\documentclass{article} %

\usepackage{amsthm}
\usepackage{graphicx}
\usepackage{enumerate}
\usepackage{mathtools}
\usepackage{xcolor} %

\usepackage{amsmath} %
\usepackage{amssymb}  %
\usepackage{bm}
\usepackage{fixme}
\fxsetup{status=draft}

\usepackage{bbm}
\usepackage{algorithm}
\usepackage{algorithmic}
\usepackage{pgfplots}
\pgfplotsset{compat=newest}

\usepackage{tikz}
\usepackage{tkz-euclide}

\usepackage{enumerate}
\usepackage{enumitem}

\usepackage{mathrsfs}
\definecolor{ao(english)}{rgb}{0.0, 0.5, 0.0}
\usepackage[colorlinks, citecolor = {ao(english)}, linkcolor = {ao(english)}]{hyperref} 
\usepackage[capitalize,nameinlink]{cleveref}
\usepackage{siunitx}
\usepackage{subcaption}

\usetikzlibrary{arrows,shapes,trees,calc,positioning,patterns,decorations.pathmorphing,decorations.markings,backgrounds}
\usetikzlibrary{matrix}

\usepgfplotslibrary{groupplots}

\newtheorem{thm}{Theorem}
\crefname{thm}{Theorem}{Theorems}

\crefname{prop}{Proposition}{Propositions}

\newtheorem{lem}{Lemma}
\crefname{lem}{Lemma}{Lemmas}

\newtheorem{cor}{Corollary}
\crefname{cor}{Corollary}{Corollaries}

\newtheorem{rem}{Remark}
\crefname{rem}{Remark}{Remark}

\newtheorem{ass}{Assumption}
\crefname{ass}{Assumption}{Assumption}

\crefname{conj}{Conjecture}{Conjectures}

\newtheorem{defn}{Definition}
\crefname{defn}{Definition}{Definitions}

\crefname{prob}{Problem}{Problems}

\crefname{appl}{Application}{Applications}

\newtheorem{exm}{Example}
\crefname{exm}{Example}{Examples}

\crefname{algorithm}{Algorithm}{Algorithms}
\crefname{paper}{Paper}{Papers}
\crefname{figure}{Figure}{Figures}
\crefname{section}{Section}{Sections}

\usepackage{dsfont}
\let\mathbb=\mathds

\newcommand{\sign}{\textnormal{sign}}
\newcommand{\supp}{\textnormal{supp}}

\newcommand{\argmin}{\operatornamewithlimits{argmin}}

\newcommand{\rel}{\textnormal{rel}}

\colorlet{FigColor1}{blue}
\colorlet{FigColor2}{red}
\colorlet{FigColor3}{ao(english)}
\colorlet{FigColor4}{orange}
\pgfplotsset{every axis plot/.append style={line width=1.5pt}}
\definecolor{bluebell}{rgb}{0.74, 0.83, 0.9}
\definecolor{airforceblue}{rgb}{0.36, 0.54, 0.66}

\crefformat{equation}{\textup{#2(#1)#3}}
\crefrangeformat{equation}{\textup{#3(#1)#4--#5(#2)#6}}
\crefmultiformat{equation}{\textup{#2(#1)#3}}{ and \textup{#2(#1)#3}}
{, \textup{#2(#1)#3}}{, and \textup{#2(#1)#3}}
\crefrangemultiformat{equation}{\textup{#3(#1)#4--#5(#2)#6}}%
{ and \textup{#3(#1)#4--#5(#2)#6}}{, \textup{#3(#1)#4--#5(#2)#6}}{, and \textup{#3(#1)#4--#5(#2)#6}}

\Crefformat{equation}{#2Equation~\textup{(#1)}#3}
\Crefrangeformat{equation}{Equations~\textup{#3(#1)#4--#5(#2)#6}}
\Crefmultiformat{equation}{Equations~\textup{#2(#1)#3}}{ and \textup{#2(#1)#3}}
{, \textup{#2(#1)#3}}{, and \textup{#2(#1)#3}}
\Crefrangemultiformat{equation}{Equations~\textup{#3(#1)#4--#5(#2)#6}}%
{ and \textup{#3(#1)#4--#5(#2)#6}}{, \textup{#3(#1)#4--#5(#2)#6}}{, and \textup{#3(#1)#4--#5(#2)#6}}

\crefdefaultlabelformat{#2\textup{#1}#3}

\usepackage{lipsum}
\usepackage{graphicx}
\usepackage{amsmath} %
\usepackage{amsthm}
\usepackage{amssymb}  %
\usepackage{fullpage}

\providecommand{\keywords}[1]{
		
	\textbf{\textit{Keywords---}} #1
}

\title{Unimodal self-oscillations and their sign-symmetry for discrete-time relay feedback systems with dead zone \thanks{The project was supported by the Israel Science Foundation (grant no. 2406/22) and the Bernard M. Gordon Center for Systems Engineering at the Technion -- IIT, while the second author was also a Jane and Larry Sherman Fellow. All the authors acknowledge the EuroTech Alliance Research Project Initiatives 2025 \& 2026 for funding this research.}}

\author{ Kang Tong\thanks{The author is with the Stephen B. Klein Faculty of Aerospace Engineering at the Technion -- Israel Institute of Technology, 3200003 Haifa, Israel. Email: {\tt\small kang.tong@technion.ac.il}} \quad Christian Grussler\thanks{The author is with the Stephen B. Klein Faculty of Aerospace Engineering at the Technion -- Israel Institute of Technology, 3200003 Haifa, Israel. Email: {\tt\small cgrussler@technion.ac.il}} \quad Michelle S. Chong \thanks{The author is with the Department of Mechanical Engineering, Eindhoven University of Technology, Eindhoven, Netherlands. Email: {\tt\small m.s.t.chong@tue.nl}}}

\begin{document}

\maketitle

\begin{abstract}
    This paper characterizes self-oscillations in discrete-time linear time-invariant (LTI) relay feedback systems with nonnegative dead zone. Specifically, we aim to establish existence criteria for unimodal self-oscillations, defined as periodic solutions where the output exhibits a single-peaked period. Assuming that the linear part of system is stable, with a strictly monotonically decreasing impulse response on its infinite support, we propose a novel analytical framework based on the theory of total positivity to address this problem. We demonstrate that unimodal self-oscillations subject to mild variation-based constraints exist only if the number of positive and negative values of the system's loop gain coincides within a given strictly positive period, i.e., the self-oscillation is sign-symmetric. Building upon these findings, we derive conditions for the existence of such self-oscillations, establish tight bounds on their periods, and address the question of their uniqueness. 
\end{abstract}


\section{Introduction}

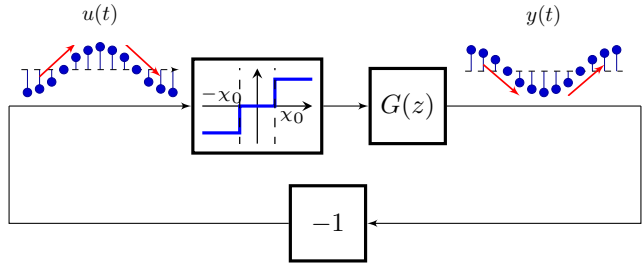
\begin{figure}[!ht]
	\centering
	\tikzstyle{neu}=[draw, very thick, align = center,circle]
	\tikzstyle{int}=[draw,minimum width=1cm, minimum height=1cm, very thick, align = center]
	\begin{tikzpicture}[>=latex',circle dotted/.style={dash pattern=on .05mm off 1.2mm,
			line cap=round}]
    \node [coordinate, name=input] {};
    \node [int, right of=input, node distance=3.3cm] (relay) {\begin{tikzpicture}[scale=1.2]
            \begin{axis}[ticks = none,width = 2.8 cm,
             ymin = -1.4,    ymax = 1.4,      %
            axis lines = middle,
            ]
                \addplot[line width = 1 pt,color = blue] table{
                    -2.5 -1
                    -0.8 -1
                    -0.8 0
                    0.8 0
                    0.8 1
                    2.5 1
                    };
    
                \draw[dashed] (-0.8,-1.5) -- (-0.8,1.5);
                \draw[dashed] (0.8,-1.5) -- (0.8,1.5);
                \node at (-1.6,0.4){\tiny$-\chi_0$};
                \node at (1.6,-0.4){\tiny$\chi_0$};
            \end{axis}
        \end{tikzpicture}};
    \node [int, right of=relay, node distance=2cm] (system) {$G(z)$};
    \node [coordinate][right of=system, node distance=3.1cm] (output) {};
    \node [int, below left of=system, node distance=1.5cm, shift={(0cm, -0.5cm)}] (feedback) {$-1$};
    \draw [-] (feedback) -| node[above] {} node[below] {} (input);
    \draw [->] (input) -- node[above] {\centering \begin{tikzpicture}[scale=0.8]
            \begin{axis}[
                width=4.2cm,
                height=2.4cm,
                axis lines=none,         %
                ticks=none,              %
                xlabel={}, ylabel={},    %
                title={$u(t)$},                %
                xlabel style={right},
                xmin=-0.5, xmax=12.5,
                ymin=-1.4, ymax=1.4,
                ]
                \addplot+[ycomb, color = blue!60!black, thin] coordinates{
                (0,-1.3)
                (1,-1.1)
                (2,-0.7)
                (3,0)
                (4,0.7)
                (5,1.1)
                (6,1.3)
                (7,1.1)
                (8,0.7)
                (9,0)
                (10,-0.7)
                (11,-1.1)
                (12,-1.3)
                (14,0)};
                \draw[dashed] (-0.5,0) -- (12.5,0);
                \draw[->, thick, red] (1,-0.4) -- (4,1.5);
                \draw[->, thick, red] (8,1.5) -- (11,-0.4);
            \end{axis}
        \end{tikzpicture}} (relay);
    \draw [->] (relay) -- node[above] {} (system);
    \draw [-] (system) -- node[above] [name=y] {\centering
        \begin{tikzpicture}[scale=0.8]
            \begin{axis}[
                width=4.2cm,
                height=2.4cm,
                axis lines=none,         %
                ticks=none,              %
                xlabel={}, ylabel={},    %
                title={$y(t)$},                %
                xlabel style={right},
                xmin=-0.5, xmax=12.5,
                ymin=-1.5, ymax=1.5,
                ]
                \addplot+[ycomb, color = blue!60!black, thin] coordinates{
                (0,1.3)
                (1,1.1)
                (2,0.7)
                (3,0)
                (4,-0.7)
                (5,-1.1)
                (6,-1.3)
                (7,-1.1)
                (8,-0.7)
                (9,0)
                (10,0.7)
                (11,1.1)
                (12,1.3)
                (14,0)};
                \draw[dashed] (-0.5,0) -- (12.5,0);
                \draw[->, thick, red] (8,-1.5) -- (11,0.4);
                \draw[->, thick, red] (1,0.4) -- (4,-1.5);
            \end{axis}
        \end{tikzpicture}} (output);
    \draw [->] (output) |- node[above] {} node[below] {} (feedback);
    \end{tikzpicture}
	
    \caption{Feedback system in discrete-time consisting of a linear time-invariant system with transfer function $G(z)$, $z\in\mathbb{C}$ and a relay with a symmetric deadzone of width $2 \chi_0$. The system here exhibits a unimodal (i.e, single-peaked) self-oscillation $u(t) = -y(t)$. \label{fig:relay_feedback}}
\end{figure}

Self-oscillatory phenomena permeate both biological and engineering domains, spanning neural networks \cite{edelstein2005limit, mackey1977oscillation}, robotic locomotion \cite{ijspeert2013dynamical}, and relay-based PID autotuning \cite{aastrom2004revisiting}.
These systems can be modelled as \emph{relay feedback systems}, where \emph{unimodality} (i.e., a single peak within a period) is a commonly analyzed oscillatory pattern \cite{tsypkin1984relay, rabi2020relay, gonccalves2001global, johansson1999fast}, as illustrated in \cref{fig:relay_feedback}.
Driven by the practical utility of relay autotuning, the prevalence of unimodality, and the intricate oscillatory dynamics inherent in relay feedback with dead zone \cite{tsypkin1984relay, rabi2020relay}, this study investigates unimodal self-oscillations in systems comprised by the feedback interconnection of a linear time-invariant system in discrete-time and a relay with a dead zone.

So far, extensive tools for the prediction of generic self-oscillations in such systems have been developed for continuous-time. This includes, on the one hand, frequency-domain graphical techniques such as describing function analysis \cite[Chap.~7.2]{khalil2002nonlinear}, Hamel locus \cite{gille1962subharmonic, le1970general}, Tsypkin locus \cite{judd1977error, tsypkin1984relay}, and the locus of a perturbed relay system (LPRS) \cite{rehan_chaotic_2024}, \cite[Chap. 2]{boiko_discontinuous_2009}, \cite[Chap. 4]{aguilar_self-oscillations_2015}, which derive \emph{approximate} conditions for the existence of self-oscillations, yet those tools fail to predict such self-oscillations according to the counterexamples from \cite{kuznetsov_harmonic_2020}.
Recent extensions of describing function analysis include replacing the sinusoidal functions with square waves \cite{chaffey_amplitude_2025}, employing the small-decay dynamic harmonic balance method \cite{zhang_transient_2025}, as well as combining it with the LPRS and a bias function \cite{rehan_analysis_2025}.

On the other hand, as self-oscillations in the \emph{output space} are, as in our setting, equivalent to limit cycles in \emph{state space}, this also includes generic state-space tools such as Poincaré maps built on the switching surface of the relay \cite{gonccalves2001global, bernardo_self-oscillations_2001, kamachkin_fixed_2022}--\cite[Chap. 3]{aguilar_self-oscillations_2015}, or the Poincaré--Bendixson theorem for planar systems \cite[Lemma 2.1]{khalil2002nonlinear}.
The Poincaré--Bendixson theorem has also been generalized to higher dimensional systems with monotone characteristics, such as monotone cyclic feedback systems \cite{mallet1990poincare}, generalized monotone systems with cones of rank $2$ \cite{sanchez_cones_2009}, 
and strongly $2$-cooperative systems \cite{weiss_generalization_2021, katz_instability_2025}, but they do not directly apply to our system due to the non-smoothness of the relay function.
Additionally, when the linear part of a relay feedback system is maximally monotone, the calculation of the limit cycle can be reformulated as the zero finding of a mixed-monotone relation \cite{das_oscillations_2022}, which can be solved by minimizing the difference of two convex functions.
Further, in case of the ideal relay, i.e., the absence of a dead zone, an existence condition, as well as an upper and lower bound for the self-oscillation with the largest period have been pursued via inspection of the open-loop step response \cite{megretski1996global}.

The analysis of self-oscillations in the discrete-time counterpart remains, to the best of our knowledge, largely unexplored. Similar to related continuous-time tools, current approaches for Lur'e feedback systems in discrete-time are either inapplicable to the relay function as they require a continuous static non-linearity \cite{paredes2024self}, or lead to conservative amplitude bounds when the dead zone is small \cite{rasvan1998self}. Note that \cite{paredes2024self, rasvan1998self} also leave the questions of unimodality and period length of a self-oscillation unaddressed, which we do so in this paper.

The purpose of this work is to develop an algebraic framework grounded in total positivity theory~\cite{karlin1968total} and its recent tractable extensions~\cite{grussler2025tractable, grussler2026dtpmp, grussler2026efficientksignconsistencyverification, grussler2026operatorsvariationboundingproperties}. This framework enables the analysis of both the \emph{existence} and the \emph{period length} of unimodal self-oscillations in discrete-time relay feedback systems without relying on differentiability assumptions or approximation techniques. Central to our approach is a rigorous classification of oscillatory patterns via \emph{cyclic variation}, defined as the number of sign changes within a single period. Unimodal self-oscillations correspond to signals whose forward differences change sign at most twice per period (see~\cref{fig:relay_feedback}). Building on this characterization, we study relay feedback systems whose open-loop gain preserves a subset of unimodal periodic signals, thereby yielding an invariant family of periodic signals that becomes amenable to a fixed-point argument. This perspective complements the variation-based state-space extensions of the Poincaré--Bendixson theorem developed for continuous-time systems in~\cite{mallet1990poincare, sanchez_cones_2009, weiss_generalization_2021, katz_instability_2025}, while shifting the focus from state-space variation to the variation of input--output sequences over one period.

We begin our investigations by studying unimodal oscillations under certain variation-based conditions and the assumption that the impulse response of the discrete-time linear time-invariant (DT-LTI) system is strictly monotonically decreasing on its infinite support. 
Based on the property of the proposed invariant set, we show that the existence of unimodal self-oscillations necessitates an equal number of positive and negative relay function output elements within one period, which is called \emph{sign-symmetry}. 
Note that a common pattern of unimodal self-oscillations in relay feedback systems, analyzed across the literatures \cite{tsypkin1984relay, gonccalves2001global, johansson1999fast, chaffey_amplitude_2025, bazanella_limit_2016, megretski1996global}, is \emph{half-wave symmetry} --- a $P$-periodic signal $u(t)$ is called half-wave symmetric when $u(t+\frac{P}{2})=-u(t)$ for all $t \in \mathbb{R}$ --- which is a special case of sign-symmetry.
Subsequently, we establish an absence condition for oscillations in DT-RFS with dead zone when the linear system has a relative degree of zero. Furthermore, we prove that the maximum oscillation period is bounded by the delay time of the LTI system. This bound serves as a prerequisite for identifying other feasible oscillation periods and can be further tightened when the impulse response is also convex on its infinite support.
The validity of these theoretical findings is illustrated through multiple numerical examples.

The present work is an extension of the authors' earlier work \cite{tong2025cdc}. The present work differs from \cite{tong2025cdc} in several major aspects: (i) The previous version only considers the special case of the relay with zero dead zone, i.e., an ideal relay. The current manuscript extends all proposed conclusions to the relay with dead zone; (ii) our earlier results have been published without proofs. These proofs can only be found in the present manuscript.

The remainder of this manuscript is organized as follows: \cref{sec:prelim} introduces the necessary preliminaries, followed by the problem formulation in \cref{sec:prob}. The characterization of self-oscillations based on variation theory is developed in \cref{sec:charac_var}. Our primary results regarding unimodal self-oscillations for plants with strictly monotonic impulse responses are derived and exemplified in \cref{sec:main_res}. Finally, \cref{sec:conclusion} concludes the paper, with all formal proofs deferred to the Appendix.

\section{Preliminaries} \label{sec:prelim}
In the following, we introduce several notations and concepts needed for our subsequent derivations and discussions. 
\subsection{Notations}
\subsubsection{Sets} We write $\mathbb{C}$ for the set of complex numbers, $\mathbb{Z}$ for the set of integers and $\mathbb{R}$ for the set of reals, with $\mathbb{Z}_{\geq0}$ ($\mathbb{Z}_{>0}$) and $\mathbb{R}_{\geq0}$ ($\mathbb{R}_{>0}$) standing for the respective subsets of nonnegative (positive) elements. For $k, l \in \mathbb{Z}$ with $k \leq l$, we write $(k:l) := \{k, k+1, \cdots, l \}$.
\subsubsection{Sequences} For a sequence $x : \mathbb{Z} \to \mathbb{R}$, with $\sup_{i \in \mathbb{Z}} |x(i)| < \infty$ or $\sum_{i \in \mathbb{Z}} |x(i)| < \infty$, we write $x \in \ell_{\infty}$ or $x \in \ell_{1}$, respectively. We use $\Delta x(t):= x(t+1)-x(t)$ to denote the \emph{forward difference operator} for a sequence $x$ and say that $x$ is \emph{convex} if $\Delta x(t) \geq \Delta x(t-1)$ for all $t \in \mathds{Z}$. Further, $x$ is called \emph{(strictly) unimodal}, if there exists an $m \in \mathbb{Z}$ for which $x(i)$ is (strictly) monotonically increasing for all $i \leq m$ and (strictly) monotonically decreasing for $t \geq m$. The \emph{support} of $x$ is defined by $\supp(x):= \{i \in \mathds{Z}: x(i) \neq 0\}$, which we say is \emph{connected} if it equals some interval $(k:l)$.

For a slice of the sequence $x$, we define the vector
\begin{align*}
    x(k:l):= \begin{bmatrix}
    x(k) & x(k+1) & \dots & x(l)
\end{bmatrix}^\top \quad \text{for } k \leq l. 
\end{align*}

\begin{figure}[!ht]
    \centering
    \begin{tikzpicture}[scale=0.9]
        \begin{axis}[height=4. cm,
            	width=9.5cm,
            	xlabel style={right},
            	xmax=12.5,xmin = 0,
                xlabel = {$i$},
                ylabel = {$x(i)$}
            	]							
            \addplot+[ycomb,black,thick]
            		coordinates{
                                    (0,0)
                                    (1,0.2)
                                    (2,0.25)
                                    (3,0.5)
            					(4,0.87)
            					(5,0.75) 
            					(6,0.65) 
            					(7,0.53) 
            					(8,0.52)
            					(9,0.47)
                                    (10,0.3)
                                    (11,0.1)
                                    (12,0.08)};
            \addplot[color = black] coordinates{(0,0) (13,0)};
            \draw[->, red] (6,0.9) -- (8,0.65)   node[right] {};
            \draw[->, red] (1,0.3) -- (3,0.75)   node[right] {};
        \end{axis}
    \end{tikzpicture}
    \caption{A single period of periodically unimodal sequence $x \in \ell_\infty(13)$, where $k_1 = 0$ and $k_2 = 4$.}
    \label{fig:PM_example}
\end{figure}
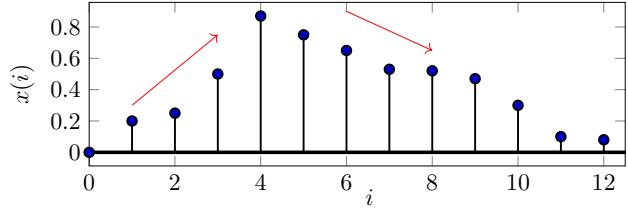

If there exists a $P \in \mathbb{Z}_{>0}$ such that $x(i) = x(i+P)$ for all $i$, then $x$ is called \emph{$P$-periodic.} The set of all bounded $P$-periodic sequences is denoted by $\ell_\infty(P)$.
Additionally, $P$-periodic $x$ is called \emph{(strictly) periodically unimodal} if there exist $k_1, k_2 \in \{0,\dots, P-1\}$ for which $x(i+k_1)$ is (strictly) monotonically increasing for all $0 \leq i \leq k_2$ and (strictly) monotonically decreasing for $k_2 \leq i \leq P-1$, as illustrated in \cref{fig:PM_example}.
Furthermore, the $P$-truncation vector of $x \in \ell_\infty(P)$ over a single period is denoted by
\begin{align*}
    x^P := \begin{bmatrix}
        x(0) \\ \vdots \\ x(P-1)
    \end{bmatrix}.
\end{align*}
For a vector $v \in \mathbb{R}^n$, we define $v_{(k)}$ to the $k$-th component by sorting all elements in $v$ in non-increasing order , i.e.,
\begin{align*}
    v_{(1)} \geq v_{(2)} \geq \dots \geq v_{(n)}.
\end{align*}
Specifically, $v_{(1)} = \max_{i\in (1:n)} v_i$ and $v_{(n)} = \min_{i\in (1:n)} v_i$.

\subsubsection{Matrices} \label{sec:preli_mat}
For matrices $M = (m_{ij}) \in \mathds{R}^{n \times m}$, we say that $M$ is
\emph{nonnegative},
$M \geq 0$ or $M \in \mathds{R}^{n \times m}_{\geq 0}$
if all elements $m_{ij} \in \mathbb{R}_{\geq 0}$ -- corresponding notations are also used for 
$M$ is \emph{nonpositive} ($M \leq 0$), \emph{positive} ($M>0$) and \emph{negative} ($M<0$) if all elements fulfill corresponding conditions. Those definitions are also used for vectors.
We use $I_n$ to denote the identity matrix in $\mathds{R}^{n \times n}$ and $\boldsymbol{c}_n$ for the vector with constant entities $c$ in $\mathbb{R}^n$.  For a vector $v = (v_i) \in \mathbb{R}^n$, the \emph{cyclic forward difference operator} $\Delta_c$ is defined as 
\begin{align*}
    \Delta_c v := \Delta_c \begin{bmatrix}
        v_1 \\ v_2 \\ \vdots \\ v_n
    \end{bmatrix} = \begin{bmatrix}
        \Delta_c v_1 \\ \Delta_c v_2 \\ \vdots \\ \Delta_c v_n
    \end{bmatrix} = \begin{bmatrix}
        v_2 - v_1 \\ v_3 - v_2 \\ \vdots \\ v_1 - v_n
    \end{bmatrix}
\end{align*}
and the \emph{circulant matrix} generated from $v$ by
\begin{align*} %
    H_v := \begin{bmatrix}
        v_1 & v_n & \cdots & v_2 \\
        v_2 & v_1 & \cdots & v_3 \\
        \vdots & \vdots &  & \vdots \\
        v_n & v_{n-1} & \cdots & v_1
    \end{bmatrix} \in \mathbb{R}^{n \times n}.
\end{align*}
Finally, we denote the \emph{cyclic back shift matrix of order $n$} by
\begin{align}
    Q_n := \begin{bmatrix} \label{eq:cyclic_back_shift_mat}
        0 & 1 \\
        I_{n-1} & 0  
    \end{bmatrix} \in \mathbb{R}^{n \times n}.
\end{align}

\subsubsection{Functions} 
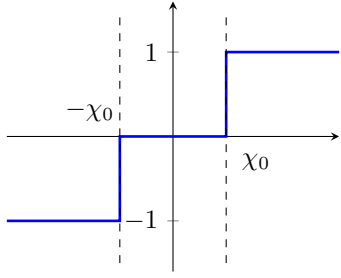
\begin{figure}[!ht]
    \centering
        \begin{tikzpicture}[>=latex']

            \begin{axis}[xtick={0},ytick={-1,1}, width =  6  cm, axis lines = middle, 
             ymin = -1.6,    ymax = 1.6,      %
            ]
                \addplot[line width = 1 pt,color = blue] table{
                    -2.5 -1
                    -0.8 -1
                    -0.8 0
                    0.8 0
                    0.8 1
                    2.5 1
                    };
    
                \draw[dashed] (-0.8,-1.5) -- (-0.8,1.5);
                \draw[dashed] (0.8,-1.5) -- (0.8,1.5);
                \node at (-1.25,0.3){$-\chi_0$};
                \node at (1.25,-0.3){$\chi_0$};
            \end{axis}
        \end{tikzpicture}
    \caption{The relay function with symmetric dead zone of width $2 \chi_0 \geq 0$.}
    \label{fig:rely_with_dead_zone}
\end{figure}

The \emph{ceiling function} $\lceil x \rceil$ maps $x \in \mathds{R}$ to the smallest integer greater than or equal to $x$.
The \emph{0-1-indicator function} $\mathds{1}_{A}$ of a set $\mathcal{A} \subset \mathbb{R}$ is defined by $\mathbb{1}_{\mathcal{A}}(t)=1$ for $t\in \mathcal{A}$ and $\mathbb{1}_{\mathcal{A}}(t)=0$ for $t \notin \mathcal{A}$.
In terms of the indicator function $\mathbb{1}_{\mathcal{A}}$, the \emph{unit pulse function} is denoted by $\delta(t):=\mathbb{1}_{\{0\}} (t)$ and the \emph{relay function with symmetric dead zone of witdh $2\chi_0 \geq 0$} in \cref{fig:rely_with_dead_zone} is denoted as 
\begin{align*}
    \rel_{\chi_0}(u) =  - \mathbb{1}_{\mathbb{R}_{< -\chi_0}}(u) + \mathbb{1}_{\mathbb{R}_{> \chi_0}}(u) , \quad u \in \mathbb{R}.
\end{align*}
If $\chi_0 = 0$, then $\sign(\cdot) := \rel_0(\cdot)$ represents the \emph{sign function} or \emph{ideal relay function}.
The number of positive, negative, and zero elements in $v \in \mathbb{R}^n$ are defined by
$P_p(v) := \sum_{i=1}^{n} \mathbb{1}_{\mathbb{R}_{> 0}} (v_i)$, $P_n(v) := \sum_{i=1}^{n} \mathbb{1}_{\mathbb{R}_{< 0}} (v_i)$ and $P_z(v) := \sum_{i=1}^{n} \mathbb{1}_{\{0\}} (v_i)$, respectively. 
Based on $P_p(\cdot)$ and $P_n(\cdot)$, we provide the precise definition of \emph{sign-symmetry}. For a vector $v \in \mathbb{R}^n$, $v$ is called sign-symmetric if $P_p(v) = P_n(v)$.

\subsection{Linear Time-Invariant Systems} 

We consider finite-dimensional {causal linear discrete-time invariant (LTI) systems} with (scalar) input $u \in \ell_\infty$ and (scalar) outputs $y \in \ell_\infty$. The \emph{impulse response} $g \in \ell_1$ is the output corresponding to the input $u(t) = \delta(t)$ and the \emph{transfer function} of the system is given by
\begin{equation} \label{eq:def_trans_fun}
G(z) = \sum_{t=0}^\infty g(t)z^{-t}, \quad z \in \mathbb{C}.
\end{equation}
We define the \emph{relative degree} $\deg(G(z))=k$ as the smallest $k\in \mathds{Z}_{\geq 0}$ such that $g(k)\neq 0$, which for rational transfer functions
\begin{equation*}
    G(z) =  \frac{\beta_m z^m+\beta_{m-1}z^{m-1}+\dots+\beta_1 z + \beta_0}{z^n+\alpha_{n-1}z^{n-1}+\dots+\alpha_1 z + \alpha_0}, \; \beta_m \neq 0, \; n \geq m
\end{equation*}
are given by $\deg(G(z)) = n-m$. The \emph{convolution operator} $\mathcal{C}_g: \ell_\infty \to \ell_\infty$ with respect to $G$ is defined by 
\begin{align}
   \mathcal{C}_g(u)(t) := (g \ast u) (t) := \sum_{\tau = -\infty}^{\infty} g(t-\tau) u(\tau), \ t \in \mathbb{Z}.
\end{align}
For $u \in \ell_\infty(P)$, it holds then that $y \in \ell_\infty(P)$ is such that
\begin{equation} \label{eq:original_conv_cal}
    \begin{aligned}
    y(t) = \mathcal{C}_g(u)(t) & = \sum_{\tau = -\infty}^{\infty} g(t-\tau) u(\tau)  \\
    & = \sum_{j=0}^{P-1} \left( u(j) \sum_{m = -\infty}^{\infty} g(t+mP-j) \right)  \\
    & = \sum_{j=0}^{P-1} u(j) \overline{g}(t-j), 
\end{aligned}
\end{equation}
with $\overline{g}(t) := \sum_{i=-\infty}^{\infty} g(t + iP)$
denoting the so-called \emph{periodic summation} of $g$. Equivalently, 
\begin{align}
    y^P =\begin{bmatrix} \label{eq:conv_to_mat_mul}
        \mathcal{C}_g(u)(0) \\ \vdots \\ \mathcal{C}_g(u)(P-1)
    \end{bmatrix} = H_{\overline{g}^P} u^P = H_{u^P} \overline{g}^P.
\end{align}

\section{Problem Statement} \label{sec:prob}
In this work, we investigate the existence of sustained unimodal oscillations in discrete-time closed-loop systems (see~\cref{fig:relay_feedback}) comprised by the feedback interconnection of a discrete linear time-invariant system $G(z)$ and relay with nonnegative dead zone. Using the introduced operator notations, this feedback interconnection is characterized by
\begin{align} \label{eq:closed_loop_sys}
    u(t) = -\mathcal{C}_g (\rel_{\chi_0}(u))(t), \quad t \in \mathbb{Z}, \, \chi_0 \geq 0.
\end{align}
Consequently, the analysis of unimodal self-oscillations is equivalent to seeking periodically unimodal fixed points $u \in \ell_\infty(P)$ satisfying \cref{eq:closed_loop_sys}.
In particular, we make the following assumption on the impulse response $g$:
\begin{ass} \label{ass:mono_dec_g}
    The impulse response $g$ to $G(z)$ satisfies $g \in \ell_{1}$, has a connected support $\supp(g)$, and is strictly monotonically decreasing on its support.  
\end{ass}
This assumption will be crucial to ensure that the system leaves the set of periodically unimodal signals invariant. Note that since $g \in \ell_1$, it also implies that $g$ is strictly positive on its support.

\section{Self-oscillations via variation} \label{sec:charac_var}
A main tool for our studies is the framework of total positivity \cite{karlin1968total,grussler2025tractable}, which characterizes convolution operators $\mathcal{C}_g$ that map periodically unimodal inputs to periodically unimodal outputs (see \cref{fig:relay_feedback}) via the notion of \emph{variation}. 
In this section, we formalize the essential definitions and review the requisite technical preliminaries associated with this framework.
\begin{defn} [Variation] \label{def:var_minus_plus}
    For a vector $v \in \mathbb{R}^n$, we define $S^-(v)$ as the number of sign changes in $v$, i.e., 
    \begin{align*}
        S^-(v) := \sum_{i=1}^{m-1} \mathbb{1}_{\mathbb{R}_{<0}} (\tilde{v}_{i} \tilde{v}_{i+1}), \quad S^-(0) = -1,
    \end{align*}
    where $\tilde{v} := [\tilde{v}_1, \tilde{v}_2, \cdots, \tilde{v}_m]^\top \in \mathbb{R}^{m}$ is the vector resulting from deleting all zeros in $v$. Further, we define $S^{+}(v) := S^{-}(\bar{v})$, where $\bar{v} \in \mathbb{R}^n$ results from replacing zero elements in $v$ by any real element that maximize $S^{-}(\bar{v})$.
\end{defn}
Note that $S^-(c) = 0$ when $c \in \mathbb{R}$ with $c \neq 0$, because $\mathbb{1}_{\mathbb{R}_{<0}} ( \tilde{v}_{2})=0$ when $\tilde{v}_{2}$ belongs to emptyset. Clearly, $S^{-}(v) \leq S^{+}(v)$, but equality does not necessarily hold. For example, if $v = [1, 0, 3]^\top$, then $S^{-}(v) = S^{-}(\tilde{v})= S^-([1,3]^\top)= 0$, but for $\bar{v} = [1, a, 3]^\top$ with $a < 0$, we obtain $S^{+}(v) = S^{-}(\bar{v}) = 2$.
We also need the notion of the so-called \emph{cyclic variation} for $v \in \mathbb{R}^n$, which is defined as follows.
\begin{align*}
    S_c^-(v) := \sup_{i \in (1:n)} S^-([v_i, \dots, v_n, v_1, \dots, v_i]), \\
    S_c^+(v) := \sup_{i \in (1:n)} S^+([v_i, \dots, v_n, v_1, \dots, v_i]).
\end{align*}
We are ready, now, to precisely define self-oscillation via  cyclic variation:
\begin{defn} [Self-oscillation] \label{def:osci_orig}
    The sequence $u \in \ell_{\infty}(P)$ is called a self-oscillation if and only if $u$ satisfies \cref{eq:closed_loop_sys} and $S_c^-(\Delta_c u^P) \geq 2$.
\end{defn}
The condition $S_c^-(\Delta_c u^P) \geq 2$ ensures that the sequence $u$ is non-constant. This is analogous to the notion of \emph{limit cycles} \cite[Chap.~2.4]{khalil2002nonlinear}, where equilibrium points, i.e., constant limit cycles, are typically. 

{Periodic unimodality of $u \in \ell_\infty(P) \setminus \{0\}$ can be equivalently using $S_c^-(\cdot)$, as illustrated in \cref{fig:S_c_example} (see \cite[Lemma~3.7]{grussler2025tractable} for the following lemma).
\begin{lem} \label{prop:PMP_three_props}
    For $u \in \ell_\infty(P) \setminus \{0\}$, the following are equivalent:
    \begin{enumerate}[label=(\roman*)]
        \item $u$ is periodically unimodal.
        \item $S_c^-(\Delta_c u^P) = 2$.
        \item $S_c^-(u^P - \gamma \boldsymbol{1}_P) \leq 2$ for all $\gamma \in \mathbb{R}$. \label{prop:PMP_three_props_case3}
    \end{enumerate}
\end{lem}

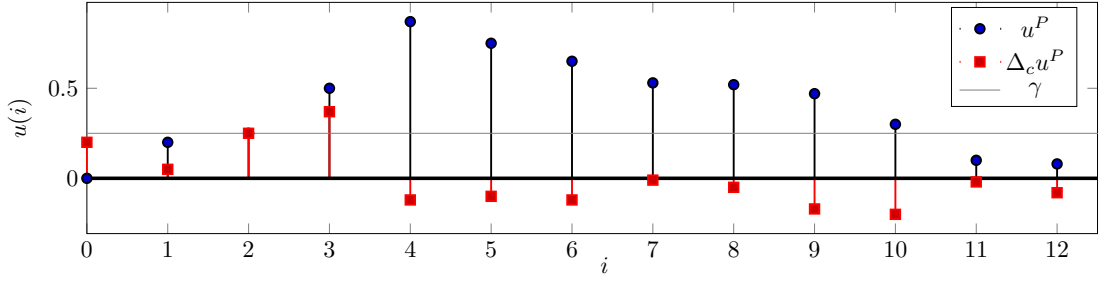
\begin{figure}[!ht]
    \centering
    \begin{tikzpicture}[scale=0.9]
        \begin{axis}[height=5 cm,
            	width= \columnwidth,
            	xlabel style={right},
            	xmax=12.5,xmin = 0,
                xlabel = {$i$},
                ylabel = {$u(i)$},
            	]							
            \addplot+[ycomb,black,thick]
            		coordinates{
                                    (0,0)
                                    (1,0.2)
                                    (2,0.25)
                                    (3,0.5)
            					(4,0.87)
            					(5,0.75) 
            					(6,0.65) 
            					(7,0.53) 
            					(8,0.52)
            					(9,0.47)
                                    (10,0.3)
                                    (11,0.1)
                                    (12,0.08)};
            \addlegendentry{$u^P$}; \label{uP}
            \addplot+[ycomb,red,thick]
            		coordinates{
                                    (0,0.2)
                                    (1,0.05)
                                    (2,0.25)
                                    (3,0.37)
            					(4,-0.12)
            					(5,-0.1) 
            					(6,-0.12) 
            					(7,-0.01) 
            					(8,-0.05)
            					(9,-0.17)
                                    (10,-0.2)
                                    (11,-0.02)
                                    (12,-0.08)}; \label{DeltauP}
            \addlegendentry{$\Delta_c u^P$};
            \addplot[color = gray, thin] coordinates{(0,0.25) (13,0.25)}; 
            \addlegendentry{$\gamma$};
            \addplot[color = black] coordinates{(0,0) (13,0)};
        \end{axis}
    \end{tikzpicture}
    \caption{Illustration of \cref{prop:PMP_three_props}: a single period of a periodically unimodal sequence $u \in \ell_\infty(13)$ and its cyclic forward difference \ref{DeltauP} $\Delta_c u^P$. Periodic unimodality is characterized by $S_c^-(\Delta u^P)=2$, or equivalently, by $S_c^-(u^P - \gamma \boldsymbol{1}_P) \leq 2$ for all $\gamma \in \mathbb{R}$.}
    \label{fig:S_c_example}
\end{figure}

In our analysis, we require the following slightly stronger form of periodic unimodality:
\begin{ass} \label{ass:osci_constr_1}
    The sequence $u \in \ell_{\infty}(P)$ satisfies 
    \begin{align} 
    S_c^-(\Delta_c u^P) = 2  \quad \text{and} \quad S_c^+(\rel_{\chi_0}(u^P)) = 2.
\end{align}
\end{ass}

The following lemma provides a way of verifying whether the circulant matrix $H_v$ satisfies $S_c^-(\Delta_c (H_v w)) \leq 2$ for any $w$ fulfilling $S_c^-(\Delta_c w) = 2$. This is important, as it allows us to check via \cref{eq:conv_to_mat_mul} whether the
\emph{loop gain} 
\begin{align} \label{eq:dt_rfs_without_delay}
    o(t) = -\mathcal{C}_g(\rel_{\chi_0}(u))(t), \quad \chi_0\geq 0,
\end{align} 
maps periodically unimodal inputs to periodically unimodal outputs --- an essential tool for our goal of finding $u$ satisfying \cref{eq:closed_loop_sys} and \cref{ass:osci_constr_1}.
\begin{lem} \label{lem:sign_PM}
    Let $v \in \mathbb{R}^n$ and $S_c^+(\rel_{\chi_0}(v)) \leq 2$ for a given $\chi_0 \geq 0$. Then, for any $w \in \mathbb{R}^n$ satisfying $S_c^-(\Delta_c w) \leq 2$, it holds that
    \begin{align*}
        S_c^-(\Delta_c H_{\rel_{\chi_0}(v)} w) \leq 2.
    \end{align*}
\end{lem}

\section{Main Results} \label{sec:main_res}
In this section, we present our main results on self-oscillations in DT-RFS with dead zone, where the linear system part fulfills \cref{ass:mono_dec_g}. We begin by deriving necessary conditions for the existence of unimodal self-oscillations satisfying \cref{ass:osci_constr_1}, where all feasible solutions $u^P$ satisfy that $\rel_{\chi_0}(u^P)$ are sign-symmetric. Based on this, we demonstrate the absence of self-oscillations for the case where $g(0) > 0$, i.e., $\deg(G(z)) = 0$. We then extend our analysis to the case $g(0) = 0$, providing lower and upper bounds for the period of self-oscillations and addressing the question of its uniqueness.

\subsection{Necessary conditions for self-oscillations}
To analyze self-oscillations in the closed-loop system, we begin by showing that if $g$ satisfies \cref{ass:mono_dec_g}, then the loop-gain \cref{eq:dt_rfs_without_delay} is invariant to signals in \cref{ass:osci_constr_1}.
\begin{thm}  \label{lem:Mickey_pres}
    Let $P > 1$, $\chi_0 \geq 0$, $g$ satisfy \cref{ass:mono_dec_g}, $u \in \ell_{\infty}(P)$ satisfy \cref{ass:osci_constr_1} and let a single period of \eqref{eq:dt_rfs_without_delay} be given by 
    \begin{align}  \label{eq:open_loop_rfs}
        o^P = -H_{\overline{g}^P}\rel_{\chi_0}(u^P).
    \end{align}
    Then, 
    \begin{align*}
        S^-_c(\Delta_c o^P) \leq 2 \text{ and } S^{-}_c(o^P) \leq 2.
    \end{align*}
    If additionally, $P_p(\rel_{\chi_0}(u^P)) = P_n(\rel_{\chi_0}(u^P))$, then $S^{+}_c(o^P) \leq 2$.
\end{thm}

The following analysis of the open-loop system properties yields our first main result: a necessary condition for the open-loop system to admit a fixed-point \cref{eq:closed_loop_sys}. When the sequence is a fixed-point for the open-loop system, the output and input sequences must share the same number of positive, negative, and zero elements, denoted by $P_p$, $P_n$, and $P_z$, in one period.

\begin{thm} \label{thm:osci_pres_open_loop}
    Let $P>1$, $g$ satisfy \cref{ass:mono_dec_g}, $u \in \ell_{\infty}(P)$ satisfy \cref{ass:osci_constr_1} and $o^P = -H_{\overline{g}^P} \rel_{\chi_0}(u^P)$ with $\chi_0 \geq 0$. 
    If $u^P$ and $o^P$ satisfy 
    \begin{align*}
        P_p(\rel_{\chi_0}(u^P)) = P_p(\rel_{\chi_0}(o^P)) \text{ and } P_n(\rel_{\chi_0}(u^P)) = P_n(\rel_{\chi_0}(o^P)), 
    \end{align*}
    then
\begin{equation*}
    P_p(\rel_{\chi_0}(u^P)) = P_n(\rel_{\chi_0}(u^P)).
\end{equation*}
If additionally, $P_z(\rel_{\chi_0}(u^P)) = 0$, then 
    \begin{align*}
        P = 2P_p(\rel_{\chi_0}(u^P)) = 2P_p(\rel_{\chi_0}(o^P)) = 2P_p(u^P) = 2P_p(o^P).
    \end{align*} 
\end{thm}

It is important to note that from \cref{lem:Mickey_pres} and \cref{thm:osci_pres_open_loop}, it can be concluded that unimodal self-oscillations satisfying \cref{ass:osci_constr_1} can only exist if the numbers of positive and negative elements of $\rel_{\chi_0}(u^P)$ in a period are equal.

\subsection{Absence of self-oscillations}
Next, we state our second main result on the absence of such unimodal self-oscillation in the case of $g(0)>0$, i.e., $\deg(G(z)) = 0$.
\begin{thm} \label{thm:non_oscillation}
    Let $g \in \ell_1$ satisfy $g(0)>0$ and \cref{ass:mono_dec_g}. Then, there does not exist any $u$ satisfying \cref{ass:osci_constr_1} such that $u = -\mathcal{C}_g(\rel_{\chi_0}(u))$ with $\chi_0 \geq 0$.
\end{thm}
An example for systems satisfying \cref{thm:non_oscillation} is the parallel interconnection of first-order lags,
\begin{align*}
    G(z) = \sum_{i=1}^{n}\frac{k_i z}{z - p_i}
\end{align*}
with $p_i \in (0, 1)$ and $k_i > 0$, where $g(0) = \sum_{i=1}^{n} k_i >0$ and the dead zone $\chi_0 = 0$. In this case, our result closely resembles the absence of self-oscillations known from continuous-time passive systems \cite[Chap.~6]{khalil2002nonlinear}.

\subsection{Self-oscillations under time-delays}

Next, we will show that self-oscillations exist if $g$ satisfies \cref{ass:mono_dec_g} with $g(0) = 0$, i.e., $\deg(G(z)) > 0$.
In this case, we can represent $G(z) = z^{-P_d} G_0(z)$, where $P_d \geq 0$ and $G_0(z)$ is causal with $g_0(0) > 0$.
This conversion is helpful for analyzing the contributions of $g_0$ and the time-delay $P_d$ independently, where the closed-loop system depicted in \cref{fig:dt_rfs_with_delay} reads as 
\begin{align} \label{eq:dt_rfs_with_delay}
    u(t) = -\mathcal{C}_g(\rel_{\chi_0}(u))(t) = -\mathcal{C}_{g_0}(\rel_{\chi_0}(u))(t-P_d),
\end{align} 
or equivalently, 
\begin{align} \label{eq:osci_delay_mat_form}
    u^P = \begin{bmatrix}
        -\mathcal{C}_{g_0}(\rel_{\chi_0}(u))(-P_d) \\ \vdots \\ -\mathcal{C}_{g_0}(\rel_{\chi_0}(u))(P-P_d-1)
    \end{bmatrix} = -Q_{P}^{P_d} H_{\overline{g}_0^P} \rel_{\chi_0}(u^P).
\end{align}
\begin{figure}[t]
	\centering
	\tikzstyle{neu}=[draw, very thick, align = center,circle]
	\tikzstyle{int}=[draw,minimum width=1cm, minimum height=1cm, very thick, align = center]
	\begin{tikzpicture}[>=latex',circle dotted/.style={dash pattern=on .05mm off 1.2mm,
			line cap=round}]
    \node [coordinate, right of=input] (sum) {};
    \node [int, right of=input, node distance=3cm] (relay) {\begin{tikzpicture}[scale=1.2]
     \begin{axis}[ticks = none,width = 2.8 cm,
             ymin = -1.4,    ymax = 1.4,      %
            axis lines = middle,
            ]
                \addplot[line width = 1 pt,color = blue] table{
                    -2.5 -1
                    -0.8 -1
                    -0.8 0
                    0.8 0
                    0.8 1
                    2.5 1
                    };
    
                \draw[dashed] (-0.8,-1.5) -- (-0.8,1.5);
                \draw[dashed] (0.8,-1.5) -- (0.8,1.5);
                \node at (-1.6,0.4){\tiny$-\chi_0$};
                \node at (1.6,-0.4){\tiny$\chi_0$};
            \end{axis}
        \end{tikzpicture}};
    \node [int, right of=relay, node distance=2cm] (system) {$G_0(z)$};
    \node [coordinate][right of=system, node distance=2cm] (output) {};
    \node [int, below of=system, node distance=1.5cm] (feedback) {$z^{-P_d}$};
    \node [int, below of=relay, node distance=1.5cm] (delay) {$-1$};
    \draw [->] (sum) -- node[above] {$u(t)$} (relay);
    \draw [->] (relay) -- node[above] {} (system);
    \draw [-] (system) -- node[above] [name=y] {$\tilde{y}(t)$} (output);
    \draw [->] (output) |- (feedback);
    \draw [->] (feedback) -- (delay);
    \draw [-] (delay) -| node[pos=0.9, left] {} node [near end, left] {} (sum);
    \end{tikzpicture}
	
    \caption{A DT Relay Feedback System with dead zone and delay module. \label{fig:dt_rfs_with_delay}}
\end{figure}
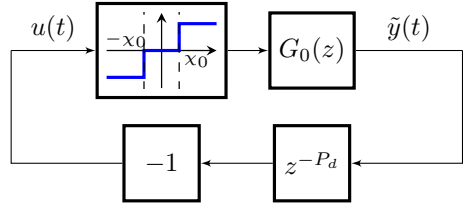
We are now ready to present our third main result: the relationship between $P_d$ and self-oscillations with period $P \geq P_d > 0$.

\begin{thm}\label{thm:bound_max_osc}
Consider system \eqref{eq:dt_rfs_with_delay} with $P_d \geq 1$, $\chi_0 \geq 0$ and $g_0$ satisfying \cref{ass:mono_dec_g}. Then, for any self-oscillation $u \in \ell_{\infty}(P)$ satisfying \cref{ass:osci_constr_1} with $P \geq P_d$, it must hold that
\begin{align*}
    2 P_d \leq P \leq 2 (P_d+P_s),
\end{align*}
where 
\begin{align*}
    P_s = \argmin_{t} \left(\sum_{k=0}^{t-1} g_0(k)- \sum_{k=t}^{\infty} g_0(k) >0 \right).
\end{align*}
In particular, $u^{2P_d}$ with $\rel_{\chi_0}(u^{2P_d}) = [\boldsymbol{1}_{P_d}^\top, -\boldsymbol{1}_{P_d}^\top]^\top$ is a self-oscillation of \eqref{eq:dt_rfs_with_delay} if and only if 
\begin{align*}
    \mathcal{C}_{g_0} (\rel_{\chi_0}{(u)}) (0) > \chi_0.
\end{align*}
\end{thm}
\begin{rem}
    A direct conclusion of \cref{thm:bound_max_osc} is that no self-oscillation satisfying \cref{ass:osci_constr_1} has a period that lies between $P_d$ and $2P_d$.
\end{rem}
Next, we show how \cref{thm:bound_max_osc} can be used to address the question of having multiple self-oscillations with possibly distinct periods.  

\begin{cor} \label{cor:max_per_to_other}
    Consider system \eqref{eq:dt_rfs_with_delay} where $P_d \geq 1$, $\chi_0 \geq 0$ and $g_0$ satisfies \cref{ass:mono_dec_g}. Further, suppose that the $P_d$-th largest element in $H_{\overline{g}^{2P_d}}  \begin{bmatrix}
            \boldsymbol{1}_{P_d} \\ -\boldsymbol{1}_{P_d}
        \end{bmatrix}$ satisfies $$\chi_0 < \left( H_{\overline{g}^{2P_d}}  \begin{bmatrix}
            \boldsymbol{1}_{P_d} \\ -\boldsymbol{1}_{P_d}
        \end{bmatrix} \right)_{(P_d)}.$$
    If there exists a self-oscillation $u^{2P_d}$ with $\rel_{\chi_0}(u^{2P_d}) = [\boldsymbol{1}_{P_d}^\top, -\boldsymbol{1}_{P_d}^\top]^\top$, then there also exist self-oscillations $u \in \ell_{\infty}(P)$ with sign pattern $\rel_{\chi_0}(u^{P}) = [\boldsymbol{1}_{\frac{P}{2}}^\top, -\boldsymbol{1}_{\frac{P}{2}}^\top]^\top$ for all (even) periods
    \begin{align*}
        P = \frac{2P_d}{2n+1} 
    \end{align*}
    whenever $P \in \mathbb{Z}_{>0}$ for $n \in \mathbb{Z}_{\geq 0}$.
\end{cor}

The following example illustrates \cref{cor:max_per_to_other} by showing how multiple self-oscillations can arise in a relay feedback system for a variety of different dead zones and time delays.
}
\begin{exm}
    Consider system \eqref{eq:dt_rfs_with_delay} and let 
    \begin{align*}
        G_0(z)= \frac{z}{z-0.1}.
    \end{align*}
    We first choose $\chi_0 = 0$ and plot the relationship between $P$ and $P_d$ by verifying \eqref{eq:osci_delay_mat_form} under the sign pattern $\sign(u^{P}) = [\boldsymbol{1}_{\frac{P}{2}}^\top, -\boldsymbol{1}_{\frac{P}{2}}^\top]^\top$ in \cref{fig:plot_osci_period} and \cref{fig:osci_period_P_d_9}. In \cref{fig:plot_osci_period}, the red points represent the existing oscillation of period $2P_d$ and the blue points represent other existing periods. 
    Under the initial conditions with different amplitudes and signs, this system displays self-oscillations with periods 2, 6, and 18 for $P_d = 9$, as illustrated in \cref{fig:osci_period_P_d_9}.

    Next, let us choose 
    \begin{align*}
        \chi_0 = 0.8 < \left(H_{\overline{g}^{6}}  \begin{bmatrix}
            \boldsymbol{1}_{3} \\ -\boldsymbol{1}_{3}
        \end{bmatrix} \right)_{(3)}\approx 0.8891 \quad \text{and} \quad P_d = 3.
    \end{align*}
    \cref{fig:osci_period_P_d_3_dead_zone} exhibits all feasible self-oscillations by setting the initial condition in our simulation, which contains all self-oscillations satisfying \cref{ass:osci_constr_1} and others. 
    We point out that the third self-oscillation in \cref{fig:osci_period_P_d_3_dead_zone} is not a part of \cref{cor:max_per_to_other}.

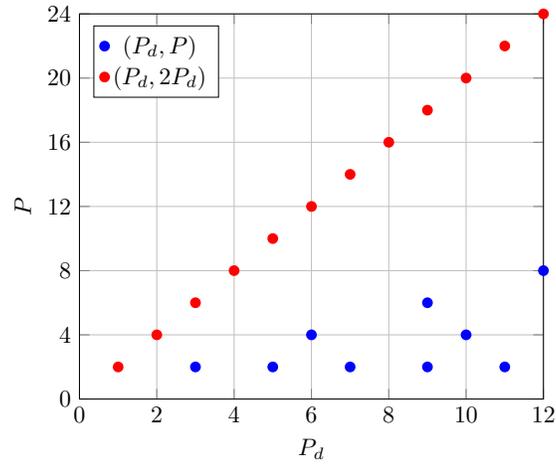
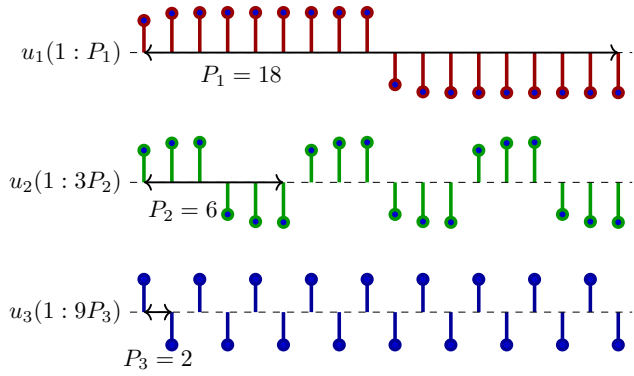
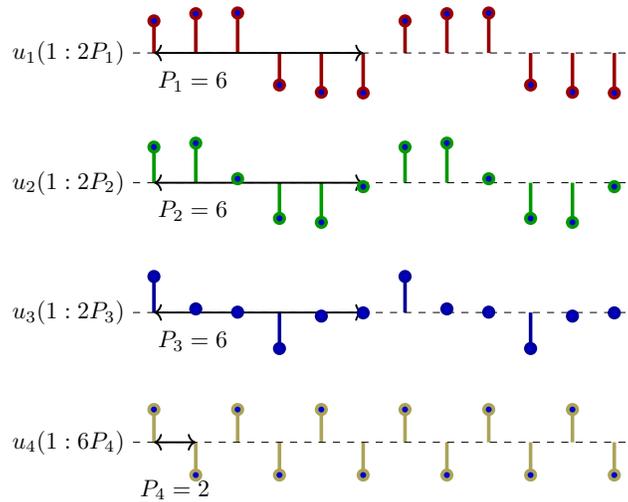
\begin{figure}[htbp]  
    \begin{subfigure}[b]{\linewidth}
        \centering
        \begin{tikzpicture}[scale=0.9]
            \begin{axis}[
                xlabel=$P_d$,
                ylabel=$P$,
                xmin=0, xmax=12,
                ymin=0, ymax=24,
                xtick={0,2,4,6,8,10,12},
                ytick={0,4,8,12,16,20,24},
                grid=both, %
                legend pos= north west, %
            ]
            \addplot[
                only marks, %
                color=blue, %
                mark=*, %
                mark size=1.5pt, %
            ] table { %
                x   y
                3   2
                5   2
                6   4
                7   2
                9   2
                9   6
                10  4
                11  2
                12  8
            };
            \addlegendentry{$(P_d,P)$}; %
    
            \addplot[
                only marks, %
                color=red, %
                mark=*, %
                mark size=1.5pt, %
            ] table { %
                x   y
                1   2
                2   4
                3   6
                4   8
                5   10
                6   12
                7   14
                8   16
                9   18
                10  20
                11  22
                12  24
            };
            \addlegendentry{$(P_d,2P_d)$}; %
            \end{axis}
        \end{tikzpicture}
        \caption{Relationship between $P$ and $P_d$ when $\chi_0 = 0$.}
        \label{fig:plot_osci_period}
    \end{subfigure}
    \par
    \vspace{1em}
    \begin{subfigure}[b]{\linewidth}
        \centering
    \begin{tikzpicture}[scale=0.9]
        \centering
        \begin{axis}[
            name=plot1,
            width=9cm, height=3.5cm,
            axis lines=none,         %
            ticks=none,              %
            xlabel={}, ylabel={},    %
            title={},                %
            axis on top=true,        %
            clip=false,              %
            xmin=0, xmax=18,
            ymin=-1.8, ymax=1.8,
        ]
        \addplot+[ycomb, color = red!60!black] 
        table {osci_P_9.txt};
        \draw[dashed] (0,0) -- (18,0);
        \node[anchor=east] at (0,0) {$u_1(1:P_1)$};
        \draw[<->, thick] (0.5,0) node[below right =2pt and 20 pt] {$P_1=18$} -- (17.5,0);
        \end{axis}

        \begin{axis}[
            at={(plot1.south)},
            name=plot2,
            anchor=north,
            width=9cm, height=3.5cm,
            axis lines=none,         %
            ticks=none,              %
            xlabel={}, ylabel={},    %
            title={},                %
            axis on top=true,        %
            clip=false,              %
            xmin=0, xmax=18,
            ymin=-1.8, ymax=1.8,
        ]
          \addplot+[ycomb, color = green!60!black] 
            table {osci_P_3.txt};
            \draw[dashed] (0,0) -- (18,0);
            \node[anchor=east] at (0,0) {$u_2(1: 3P_2)$};
            \draw[<->, thick] (0.5,0)node[below right =4pt and -2pt] {$P_2=6$} -- (5.5,0);
      \end{axis}

      \begin{axis}[
            at={(plot2.south)},
            name=plot3,
            anchor=north,
            width=9cm, height=3.5cm,
            axis lines=none,         %
            ticks=none,              %
            xlabel={}, ylabel={},    %
            title={},                %
            axis on top=true,        %
            clip=false,              %
            xmin=0, xmax=18,
            ymin=-1.8, ymax=1.8,
        ]
          \addplot+[ycomb, color = blue!60!black] 
            table {osci_P_1.txt};
            \draw[dashed] (0,0) -- (18,0);
            \node[anchor=east] at (0,0) {$u_3(1: 9P_3)$};
            \draw[<->, thick] (0.5,0) -- (1.5,0)
          node[midway, below=12pt] {$P_3=2$};
      \end{axis}    
    \end{tikzpicture}
    \caption{Three different self-oscillations when $P_d = 9$ and $\chi_0 = 0$.}
    \label{fig:osci_period_P_d_9}
    \end{subfigure}
    \par
    \vspace{1em}
    \begin{subfigure}[b]{\linewidth}
        \centering
    \begin{tikzpicture}[scale=0.9]
        \centering
        \begin{axis}[
            name=plot1,
            width=9cm, height=3.5cm,
            axis lines=none,         %
            ticks=none,              %
            xlabel={}, ylabel={},    %
            title={},                %
            axis on top=true,        %
            clip=false,              %
            xmin=0, xmax=12,
            ymin=-1.8, ymax=1.8,
        ]
        \addplot+[ycomb, color = red!60!black] 
        table {osci_dead_zone_P_3.txt};
        \draw[dashed] (0,0) -- (12,0);
        \node[anchor=east] at (0,0) {$u_1(1:2P_1)$};
        \draw[<->, thick] (0.5,0)node[below right =4pt and -2pt] {$P_1=6$} -- (5.5,0);
        \end{axis}

        \begin{axis}[
            at={(plot1.south)},
            name=plot2,
            anchor=north,
            width=9cm, height=3.5cm,
            axis lines=none,         %
            ticks=none,              %
            xlabel={}, ylabel={},    %
            title={},                %
            axis on top=true,        %
            clip=false,              %
            xmin=0, xmax=12,
            ymin=-1.8, ymax=1.8,
        ]
          \addplot+[ycomb, color = green!60!black] 
            table {osci_dead_zone_P_2_1.txt};
                
            \draw[dashed] (0,0) -- (12,0);
            \node[anchor=east] at (0,0) {$u_2(1: 2P_2)$};
            \draw[<->, thick] (0.5,0)node[below right =4pt and -2pt] {$P_2=6$} -- (5.5,0);
      \end{axis}

      \begin{axis}[
            at={(plot2.south)},
            name=plot3,
            anchor=north,
            width=9cm, height=3.5cm,
            axis lines=none,         %
            ticks=none,              %
            xlabel={}, ylabel={},    %
            title={},                %
            axis on top=true,        %
            clip=false,              %
            xmin=0, xmax=12,
            ymin=-1.8, ymax=1.8,
        ]
          \addplot+[ycomb, color = blue!60!black] 
            table {osci_dead_zone_P_2_2.txt};
                
            \draw[dashed] (0,0) -- (12,0);
            \node[anchor=east] at (0,0) {$u_3(1: 2P_3)$};
            \draw[<->, thick] (0.5,0)node[below right =4pt and -2pt] {$P_3=6$} -- (5.5,0);
      \end{axis}

      \begin{axis}[
            at={(plot3.south)},
            name=plot4,
            anchor=north,
            width=9cm, height=3.5cm,
            axis lines=none,         %
            ticks=none,              %
            xlabel={}, ylabel={},    %
            title={},                %
            axis on top=true,        %
            clip=false,              %
            xmin=0, xmax=12,
            ymin=-1.8, ymax=1.8,
        ]
          \addplot+[ycomb, color = yellow!60!black] 
            table {osci_dead_zone_P_1.txt};
            \draw[dashed] (0,0) -- (12,0);
            \node[anchor=east] at (0,0) {$u_4(1: 6P_4)$};
            \draw[<->, thick] (0.5,0) -- (1.5,0)
          node[midway, below=12pt] {$P_4=2$};
      \end{axis}

    \end{tikzpicture}
    \caption{Four different self-oscillations when $P_d = 3$ and $\chi_0 = 0.8$.}
    \label{fig:osci_period_P_d_3_dead_zone}
    \end{subfigure}
    
    \caption{Self-oscillations for Example 1.}
\end{figure}

\end{exm}

Compared with the ideal relay ($\chi_0=0$), the introduction of a non-zero dead zone $\chi_0$ may lead to the inheritance of all the self-oscillations from the ideal relay and may even introduce additional self-oscillations.

The upper bound in \cref{thm:bound_max_osc} can be estimated using time delay only, if we additionally assume that $g_0$ is convex on its support, e.g., 
\begin{align*}
    g_0(t) = a^t \text{ with } a \in (0,1) \quad \text{ for } t \in \mathbb{Z}_{\geq 0}.
\end{align*}

\begin{cor} \label{cor:convex_bound}
    Consider system \eqref{eq:dt_rfs_with_delay} where $P_d > 1$, $\chi_0 \geq 0$ and $g_0$ satisfies \cref{ass:mono_dec_g} and is convex on its support. 
    Then, for any self-oscillation $u \in \ell_{\infty}(P)$ satisfying \cref{ass:osci_constr_1} with $P \geq P_d$, it must hold that 
    \begin{align*}
        2 P_d \leq P \leq 4 P_d +2.
    \end{align*}
\end{cor}

The following example verifies that the convex impulse response function $g_0$ satisfies the tighter bound given in \cref{cor:convex_bound}.

\begin{exm}
    Consider system \eqref{eq:dt_rfs_with_delay} when $\chi_0 = 0$ with two different linear systems 
    \begin{align*}
        G_{0,1}(z)= \frac{z}{z-0.1} \text{ and } G_{0,2}(z)= \frac{z}{z-0.9},
    \end{align*} 
    with corresponding convex impulse responses
    \begin{align*}
        g_{0,1}(t) = 0.1^{t} \text{ and } g_{0,2}(t) = 0.9^{t}, \text{ for } t \in \mathbb{Z}_{\geq 0}.
    \end{align*}
    For $P \geq P_d$, we plot the relationship between $P$ and $P_d$ by verifying \eqref{eq:osci_delay_mat_form} under the sign pattern $\sign(u^{P}) = [\boldsymbol{1}_{\frac{P}{2}}^\top, -\boldsymbol{1}_{\frac{P}{2}}^\top]^\top$ in \cref{fig:plot_osci_max}. The blue points represent the largest periods of $z^{-P_d}G_{0,1}(z)$ and the red points represent the largest periods of $z^{-P_d}G_{0,2}(z)$. The two dashed lines represent the lower bound $2P_d$ and the upper bound $4P_d+2$, respectively.

\begin{figure}[htbp]
    \centering   
    \begin{tikzpicture}[scale=0.9]
        \begin{axis}[
            xlabel=$P_d$,
            ylabel=$P$,
            xmin=0, xmax=8,
            ymin=0, ymax=32,
            xtick={0,1,2,3,4,5,6,7,8},
            ytick={0,4,8,12,16,20,24,28,32},
            grid=both, %
            legend pos= north west, %
        ]
        \addplot[
            only marks, %
            color=blue, %
            mark=*, %
            mark size=1.5pt, %
        ] table { %
            x   y
            1   2
            2   4
            3   6
            4   8
            5   10
            6   12
            7   14
            8   16
        };
        \addlegendentry{$G_{0,1}(z) = \frac{z}{z-0.1}$}; %

        \addplot[
            only marks, %
            color=red, %
            mark=*, %
            mark size=1.5pt, %
        ] table { %
            x   y
            1   2
            2   6
            3   10
            4   14
            5   18
            6   22
            7   26
            8   30
        };
        \addlegendentry{$G_{0,2}(z) = \frac{z}{z-0.9}$}; %

        \addplot[magenta, dashed] coordinates{(0,2) (7.5,32)};
        \addlegendentry{$P=4P_d+2$}; %

        \addplot[cyan, dashed] coordinates{(0,0) (8,16)};
        \addlegendentry{$P=2P_d$}; %

        \end{axis}
    \end{tikzpicture}
    \caption{Relationship between $P$ and $P_d$ for different convex impulse responses under $P\geq P_d$ and $\chi_0 = 0$.}
    \label{fig:plot_osci_max}
\end{figure}
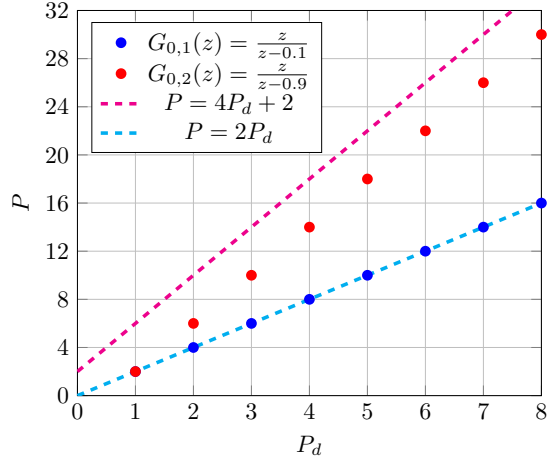

\end{exm}

\section{Conclusion} \label{sec:conclusion}
This paper has investigated the existence and characterization of unimodal self-oscillations in discrete-time relay feedback systems (DT-RFS) with a dead zone and a strictly monotonically decreasing impulse response.
By exploiting the framework of total positivity, we demonstrate that for a specific class of unimodal self-oscillations, the sign-symmetric property must hold.
Leveraging these insights, we derived conditions for the absence and existence of such oscillations as well as bounds for their periods. 
We anticipate that this discrete-time analysis will provide fundamental insights into the behavior of continuous-time systems as the limiting case of the discretization process.
Future research will focus on extending this framework to self-oscillations with higher-order variations and further elucidating the interplay between oscillation patterns, sign-symmetry, and the variation-based properties of the impulse response.

\appendix

\section{Auxiliary Results}

\begin{lem}[{\cite[Thm.~2]{grussler2025tractable}}] \label{lem:VB_2}
   Consider $v \in \mathbb{R}^n$ with $n > 3$. Then, $S_c^-(\Delta_c H_{v} w) \leq 2$ for any $w \in \mathbb{R}^n$ satisfying $S_c^-(\Delta_c w) \leq 2$, if and only if  following conditions are fulfilled:
    \begin{enumerate}
        \item $S_c^-(\Delta_c(\Delta_c v)) \leq 2$. \label{item:PMP_PM}
        \item For any $t \in (1:n)$, we have $(\Delta_c v_t)^2 \geq \Delta_c v_{t-1} \Delta_c v_{t+1}$, where the indices are taken modulo $n$; that is, $\Delta_c v_{t+kn} = \Delta_c v_{t}$ for all $k \in \mathbb{Z}$. \label{item:PMP_ineq}
    \end{enumerate}
\end{lem}
Note that if $u$ satisfies \cref{ass:osci_constr_1}, then we are only required to consider the following four main potential forms of $\rel_{\chi_0}(u^P)$: 
\begin{subequations} \label{eq:four_cases}
    \begin{align}
        & [1,\cdots,1,0,-1,\cdots,-1,0]^\top, \label{eq:case_a} \\
        & [1,\cdots,1,-1,\cdots,-1,0]^\top, \label{eq:case_b} \\
        & [1,\cdots,1,0,-1,\cdots,-1]^\top, \label{eq:case_c} \\ 
        & [1,\cdots,1,-1,\cdots,-1]^\top. \label{eq:case_d}
    \end{align}
\end{subequations}
All remaining cases can be generated by those four forms via multiplication with the cyclic back-shift matrix $Q_P^k$ for some $k \in \mathbb{Z}$. Thus, in the following proofs, we only need to discuss the four cases stated in \eqref{eq:four_cases}.

\section{Auxiliary Results}

\begin{lem}[{\cite[Thm.~2]{grussler2025tractable}}] \label{lem:VB_2}
   Consider $v \in \mathbb{R}^n$ with $n > 3$. Then, $S_c^-(\Delta_c H_{v} w) \leq 2$ for any $w \in \mathbb{R}^n$ satisfying $S_c^-(\Delta_c w) \leq 2$, if and only if  following conditions are fulfilled:
    \begin{enumerate}
        \item $S_c^-(\Delta_c(\Delta_c v)) \leq 2$. \label{item:PMP_PM}
        \item For any $t \in (1:n)$, we have $(\Delta_c v_t)^2 \geq \Delta_c v_{t-1} \Delta_c v_{t+1}$, where the indices are taken modulo $n$; that is, $\Delta_c v_{t+kn} = \Delta_c v_{t}$ for all $k \in \mathbb{Z}$. \label{item:PMP_ineq}
    \end{enumerate}
\end{lem}
Note that if $u$ satisfies \cref{ass:osci_constr_1}, then we are only required to consider the following four main potential forms of $\rel_{\chi_0}(u^P)$: 
\begin{subequations} \label{eq:four_cases}
    \begin{align}
        & [1,\cdots,1,0,-1,\cdots,-1,0]^\top, \label{eq:case_a} \\
        & [1,\cdots,1,-1,\cdots,-1,0]^\top, \label{eq:case_b} \\
        & [1,\cdots,1,0,-1,\cdots,-1]^\top, \label{eq:case_c} \\ 
        & [1,\cdots,1,-1,\cdots,-1]^\top. \label{eq:case_d}
    \end{align}
\end{subequations}
All remaining cases can be generated by those four forms via multiplication with the cyclic back-shift matrix $Q_P^k$ for some $k \in \mathbb{Z}$. Thus, in the following proofs, we only need to discuss the four cases stated in \eqref{eq:four_cases}.

\section{Proof of \cref{lem:sign_PM}}

    We first consider the cases $n=2$ and $n=3$.
    For these dimensions, every vector $v_a \in \mathbb{R}^n$ satisfies $S_c^-(v_a) \leq 2$. Hence, since $\Delta_c H_{\rel_{\chi_0}(v)} w \in \mathbb{R}^2$ or $\mathbb{R}^3$, 
    we have 
    \begin{align*}
        S_c^-(\Delta_c H_{\rel_{\chi_0}(v)} w) \leq 2 \quad \text{for any } v,w \in \mathbb{R}^2 \text{ or } \mathbb{R}^3.
    \end{align*}
Therefore, in the following, let $n>3$. Since $v \in \mathbb{R}^n$ is assumed to satisfy $S_c^+(\rel_{\chi_0}(v)) \leq 2$, we have that 
\begin{align*}
    S_c^-(\Delta_c \rel_{\chi_0}(v)) \leq 2 \quad \text{and} \quad v \neq 0.
\end{align*}
In order to prove our claim, we need to verify that $\tilde{v}:=\rel_{\chi_0}(v)$ fulfills \cref{item:PMP_PM} and \cref{item:PMP_ineq} of \cref{lem:VB_2} for the following three subcases:
    \begin{enumerate}
        \item If $S^{+}_c(\tilde{v}) = 0$, it holds that $S_c^-(\tilde{v}) =0$. Hence, $\tilde{v}$ can either be $\boldsymbol{1}_n$ or $-\boldsymbol{1}_n$. Therefore, $S_c^-(\Delta_c(\Delta_c \tilde{v})) = -1$ satisfies \cref{item:PMP_PM} of \cref{lem:VB_2} and $\Delta_c \tilde{v} = 0$ satisfies \cref{item:PMP_ineq} of \cref{lem:VB_2}.
        \item If $S^{+}_c(\tilde{v}) = 2$ and $S^{-}_c(\tilde{v}) = 0$, then after a possible cyclic permutation, $\tilde{v}$ can be either of the form $[\boldsymbol{1}_{n-1}^\top, 0]^\top$ or $[-\boldsymbol{1}_{n-1}^\top, 0]^\top$, in which case direct computations show that \cref{item:PMP_PM} and \cref{item:PMP_ineq} of \cref{lem:VB_2} are fulfilled.
        \item If $S^{+}_c(\tilde{v}) = 2$ and $S^{-}_c(\tilde{v}) = 2$ we need to check all cases given in \eqref{eq:four_cases}. We will show this only for \eqref{eq:case_a} as the others follow similarly.
        Let $\tilde{v} = [\boldsymbol{1}_{n_1}^\top, 0, -\boldsymbol{1}_{n_2}^\top, 0]^\top$ with $n_1, n_2 > 0$ and $n_1+n_2+2 = n$, then $S_c^-(\Delta_c(\Delta_c \tilde{v})) = 2$ satisfies \cref{item:PMP_PM} of \cref{lem:VB_2}. Moreover, all possible minors of order $2$ of $H_{\Delta_c \tilde{v}}$ are defined by the submatrices
        \begin{align} \label{eq:group_submat}
            \begin{bmatrix}
                0 & 0 \\ 0 & 0
            \end{bmatrix},
            \pm\begin{bmatrix}
                0 & 1 \\ 0 & 0
            \end{bmatrix},
            \pm\begin{bmatrix}
                1 & 1 \\ 0 & 1
            \end{bmatrix},
            \pm\begin{bmatrix}
                0 & 0 \\ 1 & 0
            \end{bmatrix},
            \pm\begin{bmatrix}
                1 & 0 \\ 1 & 1
            \end{bmatrix}.
        \end{align}
        As the value of $(\Delta_c (\tilde{v}_t))^2 - \Delta_c (\tilde{v}_{t-1}) \Delta_c (\tilde{v}_{t+1})$ for any $t \in (1:n)$ equals to determinat of one of those submatrices \cref{eq:group_submat}, which are all nonnegative, we have also verified \cref{item:PMP_ineq} of \cref{lem:VB_2}.
    \end{enumerate}
    Thus, in all three cases, \cref{lem:VB_2} applies, which proves our claim. \hfill $\Box$

\section{Proof of \cref{lem:Mickey_pres}}

    We begin by showing that $S_c^-(\Delta_c o^P) \leq 2$.
    Since $g$ is strictly monotonically decreasing on its support by \cref{ass:mono_dec_g}, we know that $S_c^-(\overline{g}^P)\leq S_c^+(\overline{g}^P) \leq 2$.
    Therefore, for any $u \in \ell_{\infty}(P)$ satisfying \cref{ass:osci_constr_1}, it holds that
    \begin{align*}
        S_c^-(\Delta_c o^P) & = S_c^-(-[\Delta_c \mathcal{C}_g(\rel_{\chi_0}(u))](0:P-1)) \\
        & = S_c^- \left( -\Delta_c H_{\overline{g}^P}\rel_{\chi_0}\left(u^P\right) \right)  \tag*{\text{by \cref{eq:conv_to_mat_mul}}} \\
        & = S_c^- \left( -\Delta_c H_{\rel_{\chi_0}\left(u^P\right)} \overline{g}^P \right) \leq 2. \tag*{\text{by \cref{lem:sign_PM}}}
    \end{align*}
    Further, by \cref{prop:PMP_three_props_case3} of \cref{prop:PMP_three_props} with $\gamma = 0$, it follows that $S^{-}_c(o^P)\leq 2$.

    Next, we show that $S^{+}_c(o^P) \leq 2$ under the additional assumption that $P_p(\rel_{\chi_0}(u^P)) = P_n(\rel_{\chi_0}(u^P))$. It suffices to verify this for the four cases of $\rel_{\chi_0}(u^P)$ stated in \eqref{eq:four_cases}. We will do this only for \eqref{eq:case_a} using proof by contradiction, as the others follow similarly. To this end, let $\tilde{u}^P := \rel_{\chi_0}(u^P) = [\boldsymbol{1}_{\frac{P}{2}-1}^\top, 0, -\boldsymbol{1}_{\frac{P}{2}-1}^\top, 0]^\top$ and suppose $S^+_c(o^P) > 2$. Since we have already shown that $S^-_c(o^P) \leq 2$, there exists a $t_a \in (1:P-1)$ such that $o^P_{t_a} = o^P_{t_a +1} = 0$ by the definition of $S_c^+$ in \cref{def:var_minus_plus}. Letting $\hat{g}^P = [\hat{g}^P_1, \hat{g}^P_2, \cdots, \hat{g}^P_P]^\top := (Q_P^{t_a+1})^{\top} \overline{g}^P$, it holds that
    \begin{align}
        o^P_{t_a+1} - o^P_{t_a} & = -(H_{\overline{g}^P} \tilde{u}^P)(t_a+1) + (H_{\overline{g}^P} \tilde{u}^P)(t_a) \nonumber \\
        & = -(H_{\tilde{u}^P} \overline{g}^P)(t_a+1) + (H_{\tilde{u}^P} \overline{g}^P)(t_a) \nonumber \\
        & = -(\overline{g}^P)^\top
        \left( \begin{bmatrix}
            \tilde{u}^P_{t_a+1} \\ \vdots \\ \tilde{u}^P_{t_a - P+2}
        \end{bmatrix} - 
        \begin{bmatrix}
            \tilde{u}^P_{t_a} \\ \vdots \\ \tilde{u}^P_{t_a - P+1}
        \end{bmatrix} 
         \right) \nonumber \\
        & = -(\overline{g}^P)^\top (Q_P^{t_a+1} \tilde{u}^P - Q_P^{t_a} \tilde{u}^P) \nonumber \\
        & = -(\hat{g}^P)^\top (\tilde{u}^P - Q_P^{-1}\tilde{u}^P) \nonumber \\
        & = \hat{g}^P_1 - \hat{g}^P_{\frac{P}{2}} - \hat{g}^P_{\frac{P}{2}+1} + \hat{g}^P_P, \label{eq:minus_sum}
    \end{align}
    where $\tilde{u}^P - Q_P^{-1}\tilde{u}^P =[1, \boldsymbol{0}_{\frac{P}{2}-2}^\top, -1, -1, \boldsymbol{0}_{\frac{P}{2}-2}^\top, 1]^\top$ and $Q_P$ is the cyclic back shift matrix defined in \eqref{eq:cyclic_back_shift_mat}. 
    Moreover, we also have that
    \begin{align}
         o^P_{t_a+1} + o^P_{t_a} & = -(\hat{g}^P)^\top (\tilde{u}^P + Q_P^{-1}\tilde{u}^P) \nonumber \\
         & = \hat{g}^P_1 + 2\sum_{i=2}^{\frac{P}{2}-1} \hat{g}^P_i + \hat{g}^P_{\frac{P}{2}} - \hat{g}^P_{\frac{P}{2}+1} -2\sum_{i=\frac{P}{2}+1}^{P-1} \hat{g}^P_i - \hat{g}^P_P, \label{eq:plus_sum}
    \end{align}
    where $\tilde{u}^P + Q_P^{-1}\tilde{u}^P =[1, \boldsymbol{2}_{\frac{P}{2}-2}^\top, 1, -1, -\boldsymbol{2}_{\frac{P}{2}-2}^\top, -1]^\top$. By the strict monotonicity of $g$ from \cref{ass:mono_dec_g}, there exists an $i \in (1:P)$ such that
    \begin{align*}
        \hat{g}^P_{i+1} > \dots > \hat{g}^P_P > \hat{g}^P_1 > \dots >  \hat{g}^P_{i}.
    \end{align*}
    In particular, 
    \begin{align*}
        \hat{g}^P_{\frac{P}{2}+1} > \hat{g}^P_{\frac{P}{2}} > \hat{g}^P_{P} > \hat{g}^P_{1},
    \end{align*}
    which implies that \eqref{eq:minus_sum} $< 0$, when $1 \leq i \leq \frac{P}{2}-1$
    Similarly, it holds that \eqref{eq:plus_sum} $< 0$ when $i =\frac{P}{2}$, \eqref{eq:minus_sum} $< 0$ when $\frac{P}{2}+1 \leq i \leq P-1$ and \eqref{eq:plus_sum} $> 0$ when $i =P$.
    These conditions show that \eqref{eq:minus_sum} $=0$ and \eqref{eq:plus_sum} $=0$ cannot hold simultaneously, which contradicts $o^P_{t_a}=o^P_{t_a+1}=0$.
    Hence, we conclude that $S_c^+(o^P) \leq 2$ and $S^-_c(\Delta_c o^P) \leq 2$. \hfill $\Box$

\section{Proof of \cref{thm:osci_pres_open_loop}} $P_p(\rel_{\chi_0}(u^P))=P_p(\rel_{\chi_0}(o^P))$ and $P_n(\rel_{\chi_0}(u^P))=P_n(\rel_{\chi_0}(o^P))$ hold. In order to prove that $P_p(\rel_{\chi_0}(u^P)) = P_n(\rel_{\chi_0}(u^P))$, we show that neither $P_p(\rel_{\chi_0}(u^P)) > P_n(\rel_{\chi_0}(u^P))$ nor $P_p(\rel_{\chi_0}(u^P)) < P_n(\rel_{\chi_0}(u^P))$ can arise. To this end, it suffices to consider the four cases in \eqref{eq:four_cases}. We only demonstrate this for \eqref{eq:case_a} using proof by contradiction, as the others follow similarly.
Therefore, let $\tilde{u}^P = [\boldsymbol{1}_{P_p(\tilde{u}^P)}^\top, 0, -\boldsymbol{1}_{P_n(\tilde{u}^P)}^\top, 0]^\top$, where $\tilde{u}^P :=\rel_{\chi_0}(u^P)$.
In the following, we only treat the case of $P_p(\tilde{u}^P) > P_n(\tilde{u}^P)$, because the proof of $P_p(\tilde{u}^P) < P_n(\tilde{u}^P)$ is analogous.

By \cref{lem:Mickey_pres}, we have $S_c^-(\Delta_c o^P) \leq 2$ and $S^-_c(o^P) \leq 2$, and, thus, there exists a $\tau \in \mathbb{Z}$ such that
    \begin{align} \label{eq:time_delay_mul}
        \hat{o}^P := -Q_P^{\tau} H_{\overline{g}^P} \tilde{u}^P = 
        \begin{bmatrix}
            \hat{o}^P_{+} \\ \hat{o}^P_\ominus
        \end{bmatrix},
    \end{align}
    where the circulant matrix $H_{\overline{g}^P} > 0$, $\rel_{\chi_0}(\hat{o}^P_{+}) > 0$ and $\rel_{\chi_0}(\hat{o}^P_{\ominus}) \leq 0$. 
    Note that, since the cyclic shift operater leaves $P_p(\cdot)$ and $P_n(\cdot)$ invariant, we have 
    \begin{align*}
        P_p(\hat{o}^P) = P_p(Q_P^{\tau}o^P) \quad \text{ and } \quad P_n(\hat{o}^P) = P_n(Q_P^{\tau}o^P).
    \end{align*}
    Next, by \eqref{eq:conv_to_mat_mul}, 
    \begin{align*}
        \hat{o}^P = -Q_P^{\tau} H_{\overline{g}^P} \tilde{u}^P = -Q_P^{\tau} H_{\tilde{u}^P} \overline{g}^P,
    \end{align*}
    where
    \begin{align} \label{eq:H-sign-u} 
        H_{\tilde{u}^P} = \begin{bmatrix}
            1 & 0 & -1 & \cdots & -1 & 0 & 1 & \cdots & 1 \\
            1 & 1 & 0 & \cdots & -1 & -1 & 0 & \cdots & 1 \\
            \vdots & \vdots & \vdots &  & \vdots & \vdots & \vdots &   & \vdots \\
            -1 & -1 & -1 & \cdots & 1 & 1 & 1 &\cdots & 0 \\
            0 & -1 & -1 & \cdots & 0 & 1 & 1 &\cdots & 1 \\
        \end{bmatrix}.
    \end{align}
    In particular, since $P_p(\rel_{\chi_0}(\hat{o}^P)) = P_p(\rel_{\chi_0}(o^P)) = P_p(\tilde{u}^P)$ by assumption, it follows that if $P_n(\tilde{u}^P)+1 \leq P_p(\tilde{u}^P)$ (i.e., $P_n(\tilde{u}^P) < P_p(\tilde{u}^P)$), then $\rel_{\chi_0}(\hat{o}^P_{1}) > 0$ and $\rel_{\chi_0}(\hat{o}^P_{P_n(\tilde{u}^P)+1}) > 0$.
    However, since all elements in $\overline{g}^P$ are positive by $g$ being positive on its infinite support, we arrive at our desired contradiction that
    \begin{align*}
        0 \leq \hat{o}^P_{1} + \hat{o}^P_{P_n(\tilde{u}^P)+1} =  -2\overline{g}^P_1 - \overline{g}^P_2 - \overline{g}^P_{2P_n(\tilde{u}^P)+3} -\sum_{i = 2P_n(\tilde{u}^P)+4}^P 2\overline{g}^P_i < 0.
    \end{align*}
    
    Furthermore, if $P_z(\tilde{u}^P)=0$, then $P_z(\rel_{\chi_0}(o^P))=P_z(\tilde{u}^P)=0$ and $P = 2 P_p(\tilde{u}^P) = 2 P_p(\rel_{\chi_0}(o^P)) = 2 P_p(o^P)$. \hfill $\Box$

\section{Proof of \cref{thm:non_oscillation}}

    From \cref{thm:osci_pres_open_loop}, the closed-loop system cannot admit a self-oscillation $u$ satisfying \cref{ass:osci_constr_1} if $P_p(\rel_{\chi_0}(u^P)) \neq P_n(\rel_{\chi_0}(u^P))$.
    Hence, suppose that $u^P$ satisfies \cref{ass:osci_constr_1} with $P_p(\rel_{\chi_0}(u^P)) = P_n(\rel_{\chi_0}(u^P))$ and $u^P = -H_{\overline{g}^P} \rel_{\chi_0}(u^P)$. We need to examine the four cases in \eqref{eq:four_cases}, which we will only do for \eqref{eq:case_a}, as the others follow similarly.

    To this end, we begin by noticing that if $u$  in \eqref{eq:case_a} fulfills $u=-\mathcal{C}_g((\rel_{\chi_0}u))$, then
    \begin{align*}
        u^P 
        = -H_{\overline{g}^P} \rel_{\chi_0}(u^P) = -H_{\rel_{\chi_0}(u^P)}\bar{g}^P,
    \end{align*}
    where $\rel_{\chi_0}(u^P) = [\boldsymbol{1}_{\frac{P}{2}-1}^\top, 0, -\boldsymbol{1}_{\frac{P}{2}-1}^\top, 0]^\top$ and $H_{\rel_{\chi_0}(u^P)}$ is as in \eqref{eq:H-sign-u}. Moreover, by our assumptions that $g(0)>0$ and $g$ satisfies \cref{ass:mono_dec_g}, it follows that $\overline{g}^P_1 > \overline{g}^P_2 > \dots > \overline{g}^P_P > 0$ and, therefore, 
    \begin{align*}
        0 < u^P_{\frac{P}{2}-1} = -\sum_{i = 1}^{\frac{P}{2}-1} \overline{g}^P_i + 0 \times \overline{g}^P_{\frac{P}{2}} + \sum_{i = \frac{P}{2}+1}^{P-1} \overline{g}^P_i + 0 \times \overline{g}^P_{P} < 0
    \end{align*}
   yields the desired contradiction. 
    \hfill $\Box$

\section{Proof of \cref{thm:bound_max_osc}}
    We need to verify our claim for the four main cases in \eqref{eq:four_cases}.
    Since by \cref{thm:osci_pres_open_loop}, $P_p(\rel_{\chi_0}(u^{P}))=P_n(\rel_{\chi_0}(u^{P}))$, it follows that $P$ satisfying $P \geq P_d$ must be even in case of \eqref{eq:case_a} and \eqref{eq:case_d}, and odd for \eqref{eq:case_b} and \eqref{eq:case_c}. Letting 
    \begin{equation} \label{eq:g_0_P_vector}
        \overline{g}_0^{P} = [\overline{g}^{P}_{0,1}, \dots, \overline{g}^{P}_{0,P}]^\top, 
    \end{equation}
    it follows by $g_0$ satisfying \cref{ass:mono_dec_g} and $g_0(0)>0$, that $\overline{g}^{P}_{0,i}$ is strictly monotonically decreasing in $i\in(1:P)$.
    In the following, we will consider the odd and even cases separately.  
    
    We begin with the even case, where $P = 2n \geq P_d$, $n \in \mathbb{Z}_{>0}$. In case of \eqref{eq:case_a}, we have that $\rel_{\chi_0}(u^{2n}) = [\boldsymbol{1}^{\top}_{n-1}, 0, -\boldsymbol{1}^{\top}_{n-1}, 0]^\top$ with
    \begin{align*}
        u^{2n} = -Q_{2n}^{P_d}y^{2n} \quad \text{and} \quad y^{2n} = H_{\overline{g}_0^{2n}} \rel_{\chi_0}(u^{2n}).
    \end{align*}
    Since $g_0$ satisfies \cref{ass:mono_dec_g} and $g_0(0)>0$, the vector $\overline{g}_0^{2n}$ is strictly monotonically decreasing and, therefore, $$y^{2n}_{n} = (\overline{g}^{2n}_{0, n} + \dots + \overline{g}^{2n}_{0, 2}) + 0 \times \overline{g}^{2n}_{0, 1} - (\overline{g}^{2n}_{0, 2n} + \dots + \overline{g}^{2n}_{0, n+2}) + 0 \times \overline{g}^{2n}_{0, n+1} >0.$$ Thus, since $u^{2n} = -Q_{2n}^{P_d}y^{2n}$, this implies that 
    \begin{align*}
        u^{2n}_{n+P_d} = -y^{2n}_{n} < 0.
    \end{align*}
   If $2n < n+P_d \leq n+ 2n$, it implies that
   \begin{align*}
       u^{2n}_{n+P_d} = u^{2n}_{P_d-n} <0, 
   \end{align*}
   which contradicts assumption $\rel_{\chi_0}(u^{2n}_{P_d-n}) \geq 0$ with $1<P_d-n \leq n$, which is why
   \begin{align*}
       n+P_d \leq 2n-1 \quad \text{or equivalently,} \quad n \geq P_d +1, 
   \end{align*}
   and we obtain $P = 2n \geq 2P_d+2$. 

    For \eqref{eq:case_d}, we have that $\rel_{\chi_0}(u^{2n}) = [\boldsymbol{1}^\top_{n}, -\boldsymbol{1}^\top_{n}]^\top$. Then, analogous arguments yield that 
    \begin{align*}
        u^{2n}_{n+P_d} = -y^{2n}_{n} < 0
    \end{align*}
    and, therefore, $n+P_d \leq 2n$. Hence, $P \geq 2P_d$ also in this case.
    
    Next, we consider the case of odd $P$ and since the proofs for \eqref{eq:case_b} and \eqref{eq:case_c} are similar, we restrict ourselves to \eqref{eq:case_b}. To this end, let $P = 2n+1 > P_d$ and $\rel_{\chi_0}(u^{2n+1}) = [\boldsymbol{1}_{n}^\top, 0, -\boldsymbol{1}_{n}^\top]^\top$ such that $u^{2n+1} = -Q_{2n+1}^{P_d}y^{2n+1}$ with $y^{2n+1} = H_{\overline{g}_0^{2n+1}} \rel_{\chi_0}(u^{2n+1})$. Since $P_d \geq 1$ and $g_0(0)>0$, the strict monotonicity of $\overline{g}_0^{i}$ ensures that 
    \begin{align*}
        u^{2n+1}_{n+1+P_d} = -y^{2n+1}_{n+1}=- \sum_{i=2}^{n+1}\overline{g}_0^{i} + \sum_{i=n+2}^{2n+1}\overline{g}_0^{i} < 0.
    \end{align*}
    Thus, $n+1+P_d \leq 2n+1$, i.e., $P \geq 2P_d+1$.
    
    In summary, in all cases, we have shown that $P \geq 2P_d$ and that $P = 2P_d$ can only occur in the case of \eqref{eq:case_d}.

    Next, we use proof by contradiction to show that $P \leq 2 (P_d + P_s)$. To this end, assume that $P \geq 2(P_d+P_s)+1$ with corresponding self-oscillation $u^{P} = [{u^{P}_+}^\top, {u^{P}_\ominus}^\top]^\top$, where the assumptions $ \rel_{\chi_0}(u^{P}_+)> 0$ and $\rel_{\chi_0}(u^{P}_\ominus) \leq 0$ cover all four cases in \eqref{eq:four_cases}. Then, by \cref{thm:osci_pres_open_loop}, $P_p(\rel_{\chi_0}(u^{P}))$ satisfies 
    \begin{align*}
        P_p(\rel_{\chi_0}(u^{P})) \geq \lceil \frac{P}{2} \rceil -1 \geq \lceil \frac{2(P_s+P_d)+1}{2} \rceil -1 = P_s+P_d.
    \end{align*}
    Since $u^{P} = -Q_{P}^{P_d}y^{P}$ with $y^{P} = H_{\overline{g}_0^{P}} \rel_{\chi_0}(u^{P})$, we have
    \begin{align*}
        u^{P}_{P_p(\rel_{\chi_0}(u^{P}))} = -y^{P}_{P_k} \quad \text{and} \quad \rel_{\chi_0}(u^{P}_i) =1 \quad \text{for } i \in (1:P_k),
    \end{align*}
    where 
    \begin{align*}
        P_k := P_p(\rel_{\chi_0}(u^{P})) - P_d \geq P_s,
    \end{align*} 
    it implies that
    \begin{align*}
        y^{P}_{P_k} & = \sum_{i=1}^{P_k} \rel_{\chi_0}(u^{P}_i) \overline{g}^{P}_{0,P_k+1-i} + \sum_{i=P_k+1}^{P} \rel_{\chi_0}(u^{P}_i) \overline{g}^{P}_{0,P_k+1-i} \\
        & = \sum_{i=1}^{P_k} \overline{g}^{P}_{0,P_k+1-i} + \sum_{i=P_k+1}^{P} \rel_{\chi_0}(u^{P}_i) \overline{g}^{P}_{0,P_k+1-i}.
    \end{align*} 
    Since
    \begin{align*}
        \sum_{i=1}^{P_k} \overline{g}^{P}_{0,P_k+1-i}= \sum_{i=1}^{P_k} \overline{g}^{P}_{0,i} \quad \text{and} \quad \rel_{\chi_0}(u^{P}_i) \leq 0
    \end{align*}
    for all $i \in (P_p(u^{P})+1:P)$, this implies that 
    \begin{align*}
        y^{P}_{P_k} \geq \sum_{i=1}^{P_k} \overline{g}^{P}_{0,i} - \sum_{i=P_k+1}^{P} \overline{g}^{P}_{0,i}.
    \end{align*}
    Thus, since $\overline{g}^{P}_{0,i} = \sum_{k\in \mathbb{Z}} g_0(i-1+kP)$ with $g_0(t) \geq 0$ for all $t \in \mathbb{Z}$ and $g_0 \not \equiv 0$, we have 
    \begin{align*}
        y^{P}_{P_k}> \sum_{i=0}^{P_k-1} g_0(i) - \sum_{i=P_k}^{\infty} g_0(i)  \geq \sum_{i=0}^{P_s-1} g_0(i) - \sum_{i=P_s}^{\infty} g_0(i) > 0, 
    \end{align*}
    which is a contradiction, because $-y^{P}_{P_k} = u^{P}_{P_p(u^{P})} > 0$. Hence, it must hold that $P \leq 2(P_d + P_s)$.

    Finally, we need to prove that $u^{2P_d}$ with $\rel_{\chi_0}(u^{2P_d}) = [\boldsymbol{1}_{P_d}^\top, -\boldsymbol{1}_{P_d}^\top]^\top$ is a self-oscillation of \eqref{eq:dt_rfs_with_delay} if and only if 
    \begin{align*}
        \mathcal{C}_{g_0} (\rel_{\chi_0}{(u)}) (0) > \chi_0. 
    \end{align*}
    If $(\rel_{\chi_0}{(u)})(0) > \chi_0$, then by strict monotonicity of $\overline{g}_0^{2P_d}$, we have \begin{equation*}
        \chi_0<\mathcal{C}_{g_0} (\rel_{\chi_0}{(u)})(0) \leq \dots \leq \mathcal{C}_{g_0} (\rel_{\chi_0}{(u)})(P_d-1),
    \end{equation*}and, thus, 
    \begin{align*}
        \rel_{\chi_0}(H_{\overline{g}^{2P_d}_0} \rel_{\chi_0}(u^{2P_d})) =
        \rel_{\chi_0} \left( \begin{bmatrix}
            \mathcal{C}_{g_0} (\rel_{\chi_0}{(u)})(0) \\
            \vdots \\
            \mathcal{C}_{g_0} (\rel_{\chi_0}{(u)})(2P_d-1)
        \end{bmatrix} \right) = 
        \begin{bmatrix}
            \boldsymbol{1}_{P_d} \\ -\boldsymbol{1}_{P_d}
        \end{bmatrix}.
    \end{align*}
   Consequently, 
   \begin{align*}
       \rel_{\chi_0} \left(-Q_{2P_d}^{P_d} H_{\overline{g}^{2P_d}_0} \rel_{\chi_0}(u^{2P_d})\right) = \begin{bmatrix}
            \boldsymbol{1}_{P_d} \\ -\boldsymbol{1}_{P_d}
        \end{bmatrix} = \rel_{\chi_0}(u^{2P_d})
   \end{align*}
    fulfills \eqref{eq:dt_rfs_with_delay}, i.e., $u^{2P_d}$ is a self-oscillation. Conversely, if $u^{2P_d}$ with $\rel_{\chi_0}(u^{2P_d}) = [\boldsymbol{1}_{P_d}^\top, -\boldsymbol{1}_{P_d}^\top]^\top$ is a self-oscillation of \eqref{eq:dt_rfs_with_delay}, then from 
    \begin{align*}
        u^{2P_d} = -Q_{2P_d}^{P_d} H_{\overline{g}^{2P_d}_0} \rel_{\chi_0}(u^{2P_d}), 
    \end{align*}
    it implies that $\rel_{\chi_0}(u^{2P_d}_1)=1$, i.e., $u^{2P_d}_1=\mathcal{C}_{g_0} (\rel_{\chi_0}{(u)}) (0) >\chi_0$.
    \hfill $\Box$

\section{Proof of \cref{cor:max_per_to_other}}
    We consider the self-oscillation $u^{2P_d}$ with $u^{2P_d} = -H_{\overline{g}^{2P_d}} \rel_{\chi_0}(u^{2P_d})$. Based on \cref{thm:bound_max_osc}, we know that self-oscillation $u^{2P_d}$ with $\rel_{\chi_0}(u^{2P_d}) = [\boldsymbol{1}_{P_d}^\top, -\boldsymbol{1}_{P_d}^\top]^\top$ if and only if
    \begin{align} \label{eq:C_g_0}
        \mathcal{C}_{g_0} (\rel_{\chi_0}(u))(0) & = \overline{g}^{2P_d}_{0,1} - \sum_{i=2}^{P_d+1} \overline{g}^{2P_d}_{0,i} + \sum_{i=P_d+2}^{2P_d} \overline{g}^{2P_d}_{0,i} > \chi_0.
    \end{align}
    As $P_p(\rel_{\chi_0}(u^{2P_d})) = P_n(\rel_{\chi_0}(u^{2P_d})) = P_d$, it implies that
    \begin{align*}
        u^{2P_d}_i = - u^{2P_d}_{i+P_d} \quad \text{for } i\in(1:P_d), 
    \end{align*}
    and, thus, the smallest positive element of $u^{2P_d}$ is the $P_d$-th largest element of $u^{2P_d}$, i.e., $\left( -H_{\overline{g}^{2P_d}}  \begin{bmatrix} 
            \boldsymbol{1}_{P_d} \\ -\boldsymbol{1}_{P_d}
        \end{bmatrix} \right)_{(P_d)}$, which provide a upper bound of dead zone $\chi_0$.

    Now, if $\rel_{\chi_0}(u^P) = [\boldsymbol{1}^\top_{\frac{P}{2}}, -\boldsymbol{1}^\top_{\frac{P}{2}}]^\top$ with $P = \frac{2P_d}{2n+1}$ and $n \in \mathbb{Z}_{\geq0}$, then rearranging the periodic summation yields 
    \begin{align}
        \overline{g}^{P}_{0,i} & = \sum_{j \in \mathbb{Z}} g(i-1+jP) = \sum_{j \in \mathbb{Z}} g(i-1+j\frac{2P_d}{2n+1}) \nonumber \\
        & = \sum_{l \in \mathbb{Z}} \sum_{m=0}^{2n} g\left(i - 1 + ((2n+1)l+m)\frac{2P_d}{2n+1}\right) \nonumber \\
        & = \sum_{m=0}^{2n} \sum_{l \in \mathbb{Z}} g\left(i+ \frac{2mP_d}{2n+1} -1+ 2l P_d\right) \nonumber \\
        & = \sum_{m=0}^{2n} \sum_{l \in \mathbb{Z}} g\left(i + mP -1 + 2l P_d\right) = \sum_{m=0}^{2n} \overline{g}^{2P_d}_{0,i+mP}, \label{eq:per_sum_merge}
    \end{align} 
    where $\overline{g}^{P}_{0,i}$ is defined in \cref{eq:g_0_P_vector}. Therefore,
    \begin{align*}
        \mathcal{C}_{g_0} (\rel_{\chi_0}(u^P))(0) &  = \overline{g}^{P}_{0,1} - \sum_{i=2}^{\frac{P}{2}+1} \overline{g}^{P}_{0,i} + \sum_{i=\frac{P}{2}+2}^{P} \overline{g}^{P}_{0,i} \\
        & = \sum_{k=0}^{2n} \overline{g}^{2P_d}_{0,1+kP} - \sum_{i=2}^{\frac{P}{2}+1} \overline{g}^{2P_d}_{0,i+kP}+ \sum_{i=\frac{P}{2}+2}^{P} \overline{g}^{2P_d}_{0,i+kP} \tag*{\text{by \eqref{eq:per_sum_merge}}}  \\
        &  > \overline{g}^{2P_d}_{0,1} - \sum_{i=2}^{P_d+1} \overline{g}^{2P_d}_{0,i} + \sum_{i=P_d+2}^{2P_d} \overline{g}^{2P_d}_{0,i}   \tag*{\text{by the monotonicity of $g$}}\\
        &  > \chi_0. \tag*{\text{by \cref{eq:C_g_0}}}
    \end{align*}
    Since $\overline{g}^{P}_{0,i}$ is also strictly monotonically decreasing in $i \in(1:P)$, we further have 
    \begin{align*}
        \chi_0 < \mathcal{C}_{g_0} (\rel_{\chi_0}(u))(0) \leq \mathcal{C}_{g_0} (\rel_{\chi_0}(u))(1)  \leq \dots \leq \mathcal{C}_{g_0} (\rel_{\chi_0}(u))(\frac{P}{2}-1)
    \end{align*}
    and, thus
    \begin{align*}
        \rel_{\chi_0} \left( H_{\overline{g}^{P}_0} \rel_{\chi_0}(u^P) \right) =\rel_{\chi_0} \left( \begin{bmatrix}
            \mathcal{C}_{g_0} (\rel_{\chi_0}(u))(0) \\
            \vdots \\
            \mathcal{C}_{g_0} (\rel_{\chi_0}(u))(P-1)
        \end{bmatrix} \right) = \begin{bmatrix}
            \boldsymbol{1}_{\frac{P}{2}} \\ -\boldsymbol{1}_{\frac{P}{2}}
        \end{bmatrix}.
    \end{align*}Then, by $Q_{P}^{\frac{P}{2}} = Q_{P}^{\frac{(2n+1)P}{2}} = Q_{P}^{P_d}$, we deduce that 
    \begin{align*}
        u^{P} = - Q_{P}^{\frac{P}{2}} H_{\overline{g}^{P}_0} \rel_{\chi_0}(u^P) = - Q_{P}^{P_d} H_{\overline{g}^{P}_0} \rel_{\chi_0}(u^P).
    \end{align*}
    Consequently, the closed-loop system admits a $P$-periodic self-oscillation. \hfill $\Box$

\section{Proof of \cref{cor:convex_bound}}

From \cref{thm:bound_max_osc} we already know that $P \geq 2P_d$ when $P \geq P_d$. Next, we use proof by contradiction to show that $P \leq 4P_d +2$, where want to establish the contradiction that
\begin{align*}
    u^{P}_{\lceil{\frac{P}{4}} \rceil + P_d} \leq 0 < \rel_{\chi_0} \left( u^{P}_{\lceil{\frac{P}{4}} \rceil + P_d} \right)
\end{align*}
whenever 
\begin{align*}
    P > 4 P_d + 2 \quad \text{and} \quad u^{P} = -Q_{P}^{P_d} H_{\overline{g}_0^{P}} \rel_{\chi_0}(u^{P}).
\end{align*}
It suffices to show these four cases in \eqref{eq:four_cases}, where $P$ is even in the case of \eqref{eq:case_a} and \eqref{eq:case_d}, and odd for \eqref{eq:case_b} and \eqref{eq:case_c}. 
We begin by considering even $P$ and assuming that
\begin{align} \label{eq:four_P_d_even_ass}
    P \geq 4P_d + 4 \quad \text{or equivalently,} \quad P_d \leq \frac{P-4}{4}.
\end{align}  
Then for \eqref{eq:case_d}, i.e.,  $\rel_{\chi_0}(u^{P}) = [\boldsymbol{1}^{\top}_{\frac{P}{2}},-\boldsymbol{1}^{\top}_{\frac{P}{2}}]^{\top}$, we obtain 
\begin{align*}
    u^{P}_{\lceil{\frac{P}{4}} \rceil + P_d} = -\sum_{i=1}^{\lceil{\frac{P}{4}} \rceil} \overline{g}^{P}_{0,i} + \sum_{i=\lceil{\frac{P}{4}} \rceil +1}^{\lceil{\frac{P}{4}} \rceil + \frac{P}{2}} \overline{g}^{P}_{0,i} - \sum_{i=\lceil{\frac{P}{4}} \rceil + \frac{P}{2}+1}^{P} \overline{g}^{P}_{0,i},
\end{align*} 
where $\overline{g}^{P}_{0,i}$ is defined in \cref{eq:g_0_P_vector} with $g_0$ satisfies \cref{ass:mono_dec_g}. Hence,
\begin{align*}
    u^{P}_{\lceil{\frac{P}{4}} \rceil + P_d} \leq \sum_{i=1}^{\frac{P}{2} - \lceil{\frac{P}{4}} \rceil} \left(-\overline{g}^{P}_{0, i} - \overline{g}^{P}_{0, P+1-i} + \overline{g}^{P}_{0, \frac{P}{2}+i} + \overline{g}^{P}_{0, \frac{P}{2}-i+1}\right),
\end{align*}
with equality when $\frac{P}{4} = \lceil \frac{P}{4} \rceil$; otherwise, an extra term $-\overline{g}^{P}_{0,\lceil \frac{P}{4} \rceil} + \overline{g}^{P}_{0,\lceil \frac{P}{4} \rceil+\frac{P}{2}} < 0$ appears, ensuring that the inequality is strict. Then, by the convexity assumption of $g_0$, it holds that
\begin{align*}
    \overline{g}^{P}_{0, i} + \overline{g}^{P}_{0, P+1-i} \geq \overline{g}^{P}_{0, \frac{P}{2}+i} + \overline{g}^{P}_{0, \frac{P}{2}-i+1} \quad \text{for every } i \in (1:\frac{P}{2} - \lceil{\frac{P}{4}} \rceil),
\end{align*}
which is why $u^{P}_{\lceil{\frac{P}{4}} \rceil + P_d} \leq 0$. Similarly, for the second even case \eqref{eq:case_a}, let $\rel_{\chi_0}(u^{P}) = [\boldsymbol{1}^{\top}_{\frac{P}{2}-1}, 0, -\boldsymbol{1}^{\top}_{\frac{P}{2}-1}, 0]^\top$, then
\begin{align*}
    u^{P}_{\lceil{\frac{P}{4}} \rceil + P_d} &  = -\sum_{i=1}^{\lceil{\frac{P}{4}} \rceil} \overline{g}^{P}_{0,i} + \sum_{i=\lceil{\frac{P}{4}} \rceil +1}^{\lceil{\frac{P}{4}} \rceil + \frac{P}{2}} \overline{g}^{P}_{0,i} - \sum_{i=\lceil{\frac{P}{4}} \rceil + \frac{P}{2}+1}^{P} \overline{g}^{P}_{0,i} + \overline{g}^{P}_{0,P-1} - \overline{g}^{P}_{0,\frac{P}{2}-1} \\
    & < \sum_{i=1}^{\lceil{\frac{P}{4}} \rceil} \overline{g}^{P}_{0,i} + \sum_{i=\lceil{\frac{P}{4}} \rceil +1}^{\lceil{\frac{P}{4}} \rceil + \frac{P}{2}} \overline{g}^{P}_{0,i} - \sum_{i=\lceil{\frac{P}{4}} \rceil + \frac{P}{2}+1}^{P} \overline{g}^{P}_{0,i}   \leq 0, 
\end{align*}
where the first inequality holds by $\overline{g}^{P}_{0, P-1} <\overline{g}^{P}_{0,\frac{P}{2}-1}$ and the second inequality is due to our convexity assumption on $g$.
However, for the cases \eqref{eq:case_a} and \eqref{eq:case_d}, since $P \geq 4$ and $P$ is even, it must hold that
\begin{align*}
    \lceil{\frac{P}{4}} \rceil + P_d & \leq \lceil{\frac{P}{4}} \rceil + \frac{P-4}{4}  \tag*{\text{by \cref{eq:four_P_d_even_ass}}} \\
    &\leq \frac{P+2}{4} + \frac{P-4}{4} = \frac{P-1}{2}.  \tag*{\text{due to even $P$}}
\end{align*}
Further, $\lceil{\frac{P}{4}} \rceil + P_d \leq \frac{P}{2} -1$ because $\lceil{\frac{P}{4}} \rceil + P_d$ must be an integer.
Thus, it follows that $\rel_{\chi_0} \left( u^{P}_{\lceil{\frac{P}{4}} \rceil + P_d} \right) = 1 >0$ for our given $\rel_{\chi_0}(u^P)$, which contradicts $u^{P}_{\lceil{\frac{P}{4}} \rceil + P_d} \leq 0$. This proves our claim for even $P$.

Next, we consider odd $P$, i.e., we assume that
\begin{equation*}
    P \geq 4P_d + 3 \quad \text{or equivalently,} \quad P_d \leq \frac{P-3}{4}.
\end{equation*}
Since the proofs for the odd cases \eqref{eq:case_b} and \eqref{eq:case_c} are similar, we restrict ourselves to \eqref{eq:case_b}, i.e., $\rel_{\chi_0}(u^{P}) = [\boldsymbol{1}^{\top}_{\frac{P-1}{2}}, 0, -\boldsymbol{1}^{\top}_{\frac{P-1}{2}}]^\top$, then
\begin{align*}
    u^{P}_{\lceil{\frac{P}{4}} \rceil + P_d} = -\sum_{i=1}^{\lceil{\frac{P}{4}} \rceil} \overline{g}^{P}_{0,i} + \sum_{i=\lceil{\frac{P}{4}} \rceil +1}^{\lceil{\frac{P}{4}} \rceil + \frac{P-1}{2}} \overline{g}^{P}_{0,i} - \sum_{i=\lceil{\frac{P}{4}} \rceil + \frac{P+3}{2}}^{P} \overline{g}^{P}_{0,i} \leq 0,
\end{align*} 
with equality when $\lceil \frac{P}{4} \rceil = \frac{P-1}{4}$. However, since $P_d \leq \frac{P-3}{4}$, $P \geq 3$ and $P$ is odd, it must hold that
\begin{align*}
    \lceil{\frac{P}{4}} \rceil + P_d \leq \lceil{\frac{P}{4}} \rceil + \frac{P-3}{4} \leq \frac{P+3}{4} + \frac{P-3}{4} = \frac{P}{2}.
\end{align*}
Further, $\lceil{\frac{P-1}{4}} \rceil + P_d \leq \frac{P-1}{2}$ since $\lceil{\frac{P-1}{4}} \rceil + P_d$ is an integer. Thus, it follows that $\rel_{\chi_0} \left( u^{P}_{\lceil{\frac{P}{4}} \rceil + P_d} \right) = 1 >0$ for our given $\rel_{\chi_0}(u^P)$.
Therefore, we have also shown our claim for odd $P$.
\hfill $\Box$

\bibliographystyle{elsarticle-num}
\bibliography{ref}

\end{document}